%% file: 00_main_decomdyadic.tex
\documentclass[11pt]{article}
\usepackage[utf8]{inputenc}
\usepackage[T1]{fontenc}
\usepackage{amsmath,amssymb,amsthm}
\usepackage{geometry}
\usepackage{enumerate}
\usepackage{mathtools}
\usepackage{hyperref}
\newtheorem{theorem}{Theorem}[section]
\newtheorem{corollary}{Corollary}[section]
\newtheorem{lemma}[theorem]{Lemma}
\newtheorem{example}[theorem]{Example}
\newtheorem{proposition}[theorem]{Proposition}
\theoremstyle{definition}
\newtheorem{definition}[theorem]{Definition}
\theoremstyle{remark}
\newtheorem{remark}[theorem]{Remark}

\title{Dyadic microlocal partitions for position-dependent fiber metrics and Weyl quantization}
\author{Vicente Vergara\footnote{Departamento de Matemática, Facultad de Ciencias Físicas y Matemáticas, Universidad de Concepción, Concepción, Chile. \texttt{vvergaraa@udec.cl}}}
\date{}

\begin{document}
\maketitle

\begin{abstract}
We construct a dyadic microlocal partition adapted to a position-dependent fiber metric on phase space and quantify its interaction with Weyl quantization. Under uniform ellipticity, the normalized fiber variable $\zeta=T_x\xi$ is uniformly comparable with the Euclidean frequency, so the construction does not introduce a new global symbolic order. The essential feature is instead a mesoscopic decomposition on the dyadic annulus $|\zeta|\asymp 2^k$, with block radius $R_k=2^{k/2}$. The corresponding microlocalizers satisfy explicit packing, support and derivative estimates, including the compact-local balance $|\partial_x^\alpha\partial_\xi^\beta\Lambda_{j,k}|\leq C_{\alpha,\beta,K}R_k^{|\alpha|-|\beta|}$ on each fixed spatial compact set $K$. For the spatially localized Sobolev-conjugated Weyl blocks we prove quantitative off-diagonal decay in the normalized mesoscopic distance and verify the Cotlar--Stein summability conditions, yielding strong recombination from $H^s$ to $H^{s-m_2}$. Separately, for the full product-type symbol class $S^{m_1,m_2}$ we prove the global weighted estimate $\operatorname{Op}^w(a):H^{s,\sigma}\to H^{s-m_2,\sigma-m_1}$ with an explicit finite-seminorm budget. We also record finite-order Weyl--Moyal bookkeeping and illustrate the construction through a patchwise parametrix model and the Radon transform as a model Fourier integral operator. For the localized Radon transform blocks in dimensions n=2,3, the mesoscopic construction introduces no additional dyadic localization loss beyond the standard Fourier integral operator Sobolev shift.
\end{abstract}


{\bf Keywords:} dyadic microlocal partition; position-dependent fiber metrics; mesoscopic decomposition; Weyl quantization; Calder\'on--Vaillancourt theorem; Moyal expansion; weighted Sobolev spaces; Cotlar--Stein almost orthogonality; off-diagonal decay; microlocal parametrix; Radon transform.

{\bf AMS MSC 2020:} Primary: 35S05, 35S30. Secondary: 35A18, 44A12, 81Q20, 42B20, 58J40.

\input{10_intro_dyadic}
\input{20_preli_dyadic}
\input{30_decomp_dyadic}
\input{40_quanti_dyadic}
\input{50_app1_dyadic}
\input{60_appendix_dyadic}

\bibliographystyle{plain}
\bibliography{refs}

\end{document}

%% file: 10_intro_dyadic.tex
\section{Introduction}

The purpose of this paper is to construct and quantify a dyadic microlocal partition adapted to a position-dependent fiber metric and to study its interaction with Weyl quantization.  The setting is the phase space $\mathbb{R}^n_x\times\mathbb{R}^n_\xi$, equipped with block diagonal metrics of the form
\[
g_{(x,\xi)}(y,\eta)=y^\top F y+\eta^\top G_x\eta,
\]
where $F$ is fixed, symmetric and positive definite, while $G_x$ is a smooth family of symmetric positive definite matrices satisfying uniform ellipticity and polynomial bounds on its derivatives.  We write
\[
T_x=G_x^{1/2},
\qquad
\|\xi\|_{g_x}=|T_x\xi|,
\]
and build the decomposition in the normalized fiber variable $\zeta=T_x\xi$.

A basic point, important for the scope of the paper, is that the standing ellipticity assumptions imply
\[
|T_x\xi|\asymp |\xi|
\]
uniformly in $x$.  Thus the contribution is not a new global symbolic order in the fiber variable.  The relevant geometry is instead the $x$-dependent deformation of the microlocal patches, together with the quantitative balance between the variation of $T_x$ and the spatial scale on which a frequency block can be treated uniformly.  This distinction is also essential for comparing the construction with the classical variable-scale pseudodifferential calculi discussed below.

The high-frequency decomposition is performed on the dyadic annuli
\[
\mathcal C_k=\{\zeta\in\mathbb{R}^n:2^k\leq |\zeta|<2^{k+1}\},
\]
but the individual blocks have the mesoscopic normalized radius
\[
R_k:=2^{k/2},
\qquad
R_k^2=2^k.
\]
For each $k$ we choose a maximal $R_k/2$-separated family of centers $\{\zeta_{j,k}\}_{j\in J_k}\subset\mathcal C_k$.  The resulting packing satisfies
\[
\#J_k\leq C_n R_k^n,
\]
and, after normalization by the locally finite covering sum, gives microlocalizers $\Lambda_{j,k}$ forming a partition of unity together with one low-frequency block.  On the support of a high-frequency block,
\[
|T_x\xi-\zeta_{j,k}|\leq R_k,
\qquad
\langle\xi\rangle\asymp R_k^2.
\]
The choice $R_k=2^{k/2}$ is the scale at which the two derivative directions balance: on compact spatial sets,
\[
|\partial_x^\alpha\partial_\xi^\beta\Lambda_{j,k}(x,\xi)|
\leq
C_{\alpha,\beta,K}
R_k^{|\alpha|-|\beta|}.
\]
Equivalently, one $x$-derivative costs one factor $R_k$, whereas one $\xi$-derivative gains one factor $R_k^{-1}$.  The associated reciprocal phase-space scales are encoded by
\[
g_k=R_k^2|\mathrm{d}x|^2+R_k^{-2}|\mathrm{d}\xi|^2.
\]
Globally, the same cutoff derivatives carry the polynomial factor $\langle x\rangle^{|\alpha|}$ dictated by the hypotheses on $G_x$.

This mesoscopic balance has two complementary consequences for localized symbols.  First, at the level of global finite seminorms, if $a\in S^{m_1,m_2}$ and $\imath,\ell$ are fixed nonnegative integers, then multiplication by $\Lambda_{j,k}$ is uniformly controlled after the explicit sufficient losses
\[
N_x=2\imath,
\qquad
N_\xi=\left\lceil\frac{\imath+\ell}{2}\right\rceil.
\]
Thus
\[
\|a\Lambda_{j,k}\|_{\imath,\ell}^{(m_1+N_x,m_2+N_\xi)}
\leq
C_{\imath,\ell}
\|a\|_{\imath,\ell}^{(m_1,m_2)}
\]
uniformly in $(j,k)$.  This is a finite-seminorm statement and must not be interpreted as uniform membership of all $a\Lambda_{j,k}$ in one fixed global symbol class with derivative-independent losses.  Second, on a fixed spatial compact set $K$, for $|\alpha|\leq\imath$ and $|\beta|\leq\ell$, the sharper block estimate is
\[
|\partial_x^\alpha\partial_\xi^\beta(a\Lambda_{j,k})(x,\xi)|
\leq
C_{\alpha,\beta,K}\,
\mathcal S_{\imath,\ell}(a;K)
R_k^{2m_2+|\alpha|-|\beta|}.
\]
This latter estimate is the one that governs the natural Sobolev mapping and the almost-orthogonality argument.

After Weyl quantization, the spatially localized blocks are conjugated to a common $L^2$ space.  For fixed compactly supported cutoffs $\psi,\psi'$ we set
\[
T_{j,k}=\psi\operatorname{Op}^w(a\Lambda_{j,k})\psi',
\qquad
\widetilde T_{j,k}
=
\langle D\rangle^{s-m_2}T_{j,k}\langle D\rangle^{-s}.
\]
The mesoscopic rescaling gives a uniform diagonal bound for $\widetilde T_{j,k}$.  More importantly, if $\lambda=(j,k)$ and $\mu=(j',k')$, with $R_\lambda=R_k$ and $\zeta_\lambda=\zeta_{j,k}$, the natural normalized block distance is
\[
\operatorname{dist}_{\mathrm{mes}}(\lambda,\mu)
:=
\frac{|\zeta_\lambda-\zeta_\mu|}{R_\lambda+R_\mu}.
\]
The kernels are localized at spatial scale $R_\lambda^{-1}$, and a nonstationary-phase argument yields, for every $M\in\mathbb{N}_0$,
\[
\|\widetilde T_\lambda^*\widetilde T_\mu\|_{L^2\to L^2}
+
\|\widetilde T_\lambda\widetilde T_\mu^*\|_{L^2\to L^2}
\leq
C_M
\bigl(1+\operatorname{dist}_{\mathrm{mes}}(\lambda,\mu)\bigr)^{-M}.
\]
The packing of the mesoscopic nets implies summability for powers larger than $n$.  Choosing
\[
M_{\operatorname{CS}}=2n+2
\]
therefore verifies the square-root summability required by Cotlar--Stein.  Consequently, the high-frequency family recombines strongly from $H^s$ to $H^{s-m_2}$; the recombination theorem is quantitative rather than conditional.  The appendix records an explicit sufficient finite derivative budget for this argument.

The paper also separates this blockwise mesoscopic mechanism from a genuinely global mapping property of the full symbol class.  For
\[
\|u\|_{H^{s,\sigma}}
:=
\|\langle x\rangle^\sigma\langle D\rangle^s u\|_{L^2},
\]
the product-type class $S^{m_1,m_2}$ satisfies
\[
\operatorname{Op}^w(a):
H^{s,\sigma}(\mathbb{R}^n)
\longrightarrow
H^{s-m_2,\sigma-m_1}(\mathbb{R}^n).
\]
This estimate is global in $x$ and requires no spatial cutoff and no additional decay hypothesis on $G_x$.  Its proof uses the exact Weyl symbol of the conjugated operator, for which the spatial and frequency orders cancel, followed by a global Calder\'on--Vaillancourt estimate.  The dependence on $a$ is again through an explicit finite seminorm budget.  This global statement concerns the complete symbol $a$; under the present assumptions it is not asserted uniformly at the same natural orders for each individual factor $a\Lambda_{j,k}$, whose global $x$-derivatives inherit polynomial growth from the moving normalization.

The geometric background of this construction belongs to a classical line of pseudodifferential analysis.  Beals and Fefferman developed calculi for spatially inhomogeneous pseudodifferential operators in which the admissible local scales vary over phase space \cite{BealsFefferman1973,BealsFefferman1974,BealsGeneralCalculus1975}.  Hörmander's Weyl calculus subsequently placed variable phase-space metrics at the center of a broad symbolic framework \cite{HormanderWeyl1979,HormanderIII}.  The present paper does not propose a replacement or a generalization of those calculi.  In particular, because $|T_x\xi|\asymp|\xi|$, the moving fiber normalization does not create a new global order structure.  Instead, we specialize to the concrete deformation $\zeta=T_x\xi$ and extract from it an explicit dyadic and mesoscopic partition, with quantitative packing, derivative balance, kernel localization, block distance and off-diagonal attenuation.  At the level of scaling, the metric $g_k=R_k^2|\mathrm{d}x|^2+R_k^{-2}|\mathrm{d}\xi|^2$ expresses reciprocal spatial and frequency resolutions; our contribution is to realize this balance constructively for the moving fiber geometry and to make the resulting recombination estimates explicit.

This positioning also clarifies the relation with more standard ingredients of the paper.  Weyl composition and the Moyal product are used in their usual form, and the $L^2$ estimates rely on Calder\'on--Vaillancourt; see, for example, \cite{HormanderIII,TaylorPDEII,FollandHAPS}.  The scale-adapted rescaling is compatible with standard semiclassical reasoning \cite{ZworskiSemiclassical}.  Cotlar--Stein supplies the abstract almost-orthogonality principle, but the required interaction decay is proved here from the mesoscopic kernel geometry rather than imposed as a separate hypothesis.  The novelty is therefore not a new general symbolic calculus, but the explicit quantitative bridge from the moving normalization $T_x$ to a dyadic phase-space partition and a verifiable operator recombination mechanism.

The use of scale-dependent decompositions is consistent with other parts of harmonic and microlocal analysis where adapted coverings simplify local estimates and recombination.  Multiparameter decompositions, directional systems such as curvelets, and parabolic decompositions for Fourier integral operators all exploit this principle; see, for instance, \cite{TanakaYabuta2019,Wang2021,CandesDonohoCurvelets,GuoLabateFIO,Rozendaal2021,LiuRozendaalSongYan2024}.  The construction here differs in that the frequency deformation is generated by the base-point-dependent map $T_x$ and the block radius is tied to the dyadic level through the mesoscopic relation $R_k^2=2^k$.

The applications at the end of the paper should be read in this restricted sense.  They illustrate how the partition can be inserted into standard microlocal arguments, but they do not claim that the uniformly elliptic regime produces symbolic orders beyond the usual Euclidean classes.  The parametrix discussion is a patchwise model based on local symbolic inversion, exact Weyl composition modulo a smoothing remainder, and strong mesoscopic recombination.  The Radon transform is used as a model Fourier integral operator to test the local--global bookkeeping of the construction rather than as a replacement for the classical Radon theory.  For the localized Radon blocks the compensated estimate produces no additional mesoscopic localization loss, namely $N_{\operatorname{loc}}=0$.  Thus, with $m_R=(n-1)/2$, the low-dimensional cases read
\[
n=2:\quad H^s_x\to H^{s+\frac12-m_2}_{s,\omega},
\qquad
n=3:\quad H^s_x\to H^{s+1-m_2}_{s,\omega}.
\]

\paragraph{Organization of the paper.}
Section~\ref{sec:metricas} fixes the metric assumptions and the product-type symbol classes.  It introduces $T_x=G_x^{1/2}$, proves the required smoothness estimates for $T_x$, and records the uniform comparison $|T_x\xi|\asymp|\xi|$.

Section~\ref{subsec:dyadic-decomposition} constructs the mesoscopic dyadic partition.  The normalized annuli have radius of order $2^k=R_k^2$, while the individual blocks have radius $R_k=2^{k/2}$.  The section proves the packing bound $\#J_k\leq C_nR_k^n$, the exact partition of unity, the support localization, and the derivative estimate of Proposition~\ref{prop:Lambda-deriv-bounds:full}.

Section~\ref{sec:quantization-local-global} contains the quantitative operator theory.  Proposition~\ref{prop:mult-by-Lambda-preserves-S-anisotropic:full} gives the global finite-seminorm localization bound, while the local estimates exploit the sharper mesoscopic derivative balance.  The Sobolev-conjugated blocks are localized at spatial scale $R_k^{-1}$, and Lemma~\ref{lem:mesoscopic-offdiagonal-cotlar} proves rapid decay in the normalized block distance together with the discrete summability needed for Cotlar--Stein.  Theorem~\ref{thm:global-weighted-sobolev} gives the global weighted estimate
\[
\operatorname{Op}^w(a):H^{s,\sigma}\to H^{s-m_2,\sigma-m_1},
\]
and Theorem~\ref{thm:global-sum-semiclassical} gives the quantitative recombination of the spatially localized mesoscopic blocks.

Section~\ref{sec:applications} contains two model uses of the localized construction.  Section~\ref{sec:parametrix} presents a patchwise parametrix model based on local symbolic inversion, exact Weyl composition modulo smoothing, and the strong mesoscopic recombination theorem.  Section~\ref{sec:radon-dyadic} treats the Radon transform as a model Fourier integral operator, proves the compensated local bound with $N_{\operatorname{loc}}=0$ for $n=2,3$, and obtains the global Radon recombination by composition with the strongly convergent pseudodifferential block sum.

Appendix~\ref{sec:appendix} collects Weyl--Moyal composition, Calder\'on--Vaillancourt, Cotlar--Stein, and the finite derivative bookkeeping used in the global weighted Sobolev and mesoscopic almost-orthogonality arguments.

%% file: 20_preli_dyadic.tex
\section{Metric assumptions and preliminary estimates}\label{sec:metricas}

This section fixes the notation, the metric hypotheses, and the elementary estimates used throughout the paper.  The main point is that the fiber metric is uniformly equivalent to the Euclidean one, while the normalization map $T_x=G_x^{1/2}$ may vary with $x$.  The derivative estimates for $T_x$ are the source of the finite losses appearing later in the dyadic localization.

\subsection*{Notation and conventions}

We use standard multi-index notation.  If $\gamma=(\gamma_1,\dots,\gamma_n)\in\mathbb{N}_0^n$, then
\[
\gamma!:=\prod_{j=1}^n\gamma_j!,
\qquad
\binom{\gamma}{\alpha}:=\frac{\gamma!}{\alpha!(\gamma-\alpha)!}.
\]
The Leibniz rule is written as
\[
\partial^\gamma(fg)
=
\sum_{\alpha\leq\gamma}
\binom{\gamma}{\alpha}
\partial^\alpha f\,\partial^{\gamma-\alpha}g.
\]

We write $|\cdot|$ for the Euclidean norm and
\[
\langle x\rangle=(1+|x|^2)^{1/2},
\qquad
\langle \xi\rangle=(1+|\xi|^2)^{1/2}.
\]
For a linear map $A$, $\|A\|_{\mathrm{op}}$ denotes its Euclidean operator norm.  More generally, if $B$ is a multilinear map, then $\|B\|_{\mathrm{op}}$ denotes the supremum of $|B(v_1,\dots,v_k)|$ over $|v_j|\leq 1$.

All constants denoted by $C$ may change from line to line.  Constants with subscripts, such as $C_{\alpha,\beta}$, may depend on the displayed indices and on the structural constants in the metric assumptions, but never on $x$, $\xi$ or the dyadic parameters introduced later.

\begin{definition}[Product symbol classes]\label{def:anisotropic-symbols}
Let $m_1,m_2\in\mathbb{R}$.  The class $S^{m_1,m_2}$ consists of all functions $a\in C^\infty(\mathbb{R}^n_x\times\mathbb{R}^n_\xi)$ such that, for every pair of multi-indices $\alpha,\beta$, there exists $C_{\alpha,\beta}>0$ with
\[
|\partial_x^\alpha\partial_\xi^\beta a(x,\xi)|
\le
C_{\alpha,\beta}
\langle x\rangle^{m_1-|\alpha|}
\langle \xi\rangle^{m_2-|\beta|},
\qquad
(x,\xi)\in\mathbb{R}^{2n}.
\]
For nonnegative integers $\imath,\ell$, the corresponding truncated seminorm is
\[
\|a\|_{\imath,\ell}^{(m_1,m_2)}
:=
\sup_{\substack{|\alpha|\leq \imath\\|\beta|\leq \ell}}
\sup_{(x,\xi)\in\mathbb{R}^{2n}}
\frac{
|\partial_x^\alpha\partial_\xi^\beta a(x,\xi)|
}{
\langle x\rangle^{m_1-|\alpha|}
\langle \xi\rangle^{m_2-|\beta|}
}.
\]
\end{definition}

\begin{remark}
The terminology ``anisotropic'' in the sequel refers to the $x$-dependent deformation of the fiber variable through $T_x\xi$, not to a new global symbolic order.  Under the uniform ellipticity assumptions below, $\|T_x\xi\|$ is uniformly comparable with $|\xi|$, and hence the symbolic weights are equivalent to the Euclidean ones.
\end{remark}

\subsection{Metric hypotheses}

Let
\[
G\in C^\infty(\mathbb{R}^n;\operatorname{Sym}^+(n))
\]
be a smooth family of symmetric positive definite matrices.  We assume that there are constants
\[
0<\lambda_G^{\min}\leq \lambda_G^{\max}<\infty
\]
such that
\[
\sigma(G_x)\subset [\lambda_G^{\min},\lambda_G^{\max}]
\qquad
\text{for every } x\in\mathbb{R}^n.
\]
We also assume that, for every multi-index $\beta$, there exists $c_\beta>0$ such that
\begin{equation}\label{eq:bd-G-beta}
\|\partial_x^\beta G_x\|_{\mathrm{op}}
\leq
c_\beta \langle x\rangle^{|\beta|},
\qquad
x\in\mathbb{R}^n.
\end{equation}

We denote by
\[
T_x:=G_x^{1/2}
\]
the principal positive square root of $G_x$.  The fiber quadratic form and its associated norm are
\[
g_x(\xi):=\xi^\top G_x\xi=|T_x\xi|^2,
\qquad
\|\xi\|_{g_x}:=|T_x\xi|.
\]

We also fix a constant matrix $F\in\operatorname{Sym}^+(n)$ and consider the phase-space metric
\[
g_{(x,\xi)}(y,\eta)
=
y^\top F y+\eta^\top G_x\eta.
\]
Its inverse quadratic form is
\[
g^{(x,\xi)}(p,q)
=
p^\top F^{-1}p+q^\top G_x^{-1}q.
\]

\begin{lemma}[Derivatives of the matrix square root]\label{lem:matrix-sqrt-deriv}
For every multi-index $\gamma$, there exists a constant $C_\gamma>0$, depending only on $\lambda_G^{\min}$, $\lambda_G^{\max}$, $|\gamma|$, and the constants $c_\beta$ with $|\beta|\leq |\gamma|$, such that
\[
\|\partial_x^\gamma T_x\|_{\mathrm{op}}
\leq
C_\gamma\langle x\rangle^{|\gamma|},
\qquad
x\in\mathbb{R}^n.
\]
For $\gamma=0$, one may take $C_0=(\lambda_G^{\max})^{1/2}$.
\end{lemma}

\begin{proof}
The case $\gamma=0$ follows from the spectral bound on $G_x$.  Assume now that $|\gamma|\geq 1$.

Choose a closed contour $\Gamma\subset\mathbb{C}\setminus(-\infty,0]$ enclosing the interval $[\lambda_G^{\min},\lambda_G^{\max}]$ and staying at positive distance from it.  Since $\sigma(G_x)\subset[\lambda_G^{\min},\lambda_G^{\max}]$ for all $x$, there exists $\delta>0$ such that
\[
\operatorname{dist}(\Gamma,\sigma(G_x))\geq \delta
\qquad
\text{for all } x\in\mathbb{R}^n.
\]
The principal square root has the Dunford representation
\[
T_x
=
G_x^{1/2}
=
\frac{1}{2\pi i}
\int_\Gamma
\lambda^{1/2}(\lambda I-G_x)^{-1}\,\mathrm{d}\lambda.
\]
For $\lambda\in\Gamma$,
\[
\|(\lambda I-G_x)^{-1}\|_{\mathrm{op}}
\leq
\delta^{-1}.
\]

Differentiating the resolvent identity iteratively, $\partial_x^\gamma(\lambda I-G_x)^{-1}$ is a finite sum of terms of the form
\[
(\lambda I-G_x)^{-1}
(\partial_x^{\alpha_1}G_x)
(\lambda I-G_x)^{-1}
\cdots
(\partial_x^{\alpha_p}G_x)
(\lambda I-G_x)^{-1},
\]
where $1\leq p\leq |\gamma|$, $\alpha_1+\cdots+\alpha_p=\gamma$, and $|\alpha_j|\geq 1$.  Hence
\[
\|\partial_x^\gamma(\lambda I-G_x)^{-1}\|_{\mathrm{op}}
\leq
C_\gamma'
\langle x\rangle^{|\gamma|},
\]
where $C_\gamma'$ depends only on $\delta$, $|\gamma|$, and the constants $c_\beta$ with $|\beta|\leq|\gamma|$.  Inserting this bound in the Dunford formula gives
\[
\|\partial_x^\gamma T_x\|_{\mathrm{op}}
\leq
\frac{1}{2\pi}
\operatorname{length}(\Gamma)
\sup_{\lambda\in\Gamma}|\lambda^{1/2}|
C_\gamma'
\langle x\rangle^{|\gamma|}.
\]
This proves the claim.
\end{proof}

\begin{proposition}[Algebraic comparison of the metric]\label{prop:gx-basic}
Under the spectral assumptions on $G_x$, the following properties hold.

\begin{enumerate}
\item For each $x$, $\|\cdot\|_{g_x}$ is a norm on $\mathbb{R}^n$.
\item For every $x,\xi$,
\[
(\lambda_G^{\min})^{1/2}|\xi|
\leq
\|\xi\|_{g_x}
\leq
(\lambda_G^{\max})^{1/2}|\xi|.
\]
\item For every $x,y,\xi$,
\[
\left(\frac{\lambda_G^{\min}}{\lambda_G^{\max}}\right)^{1/2}
\|\xi\|_{g_y}
\leq
\|\xi\|_{g_x}
\leq
\left(\frac{\lambda_G^{\max}}{\lambda_G^{\min}}\right)^{1/2}
\|\xi\|_{g_y}.
\]
\item If $\lambda_F^{\min}$ and $\lambda_F^{\max}$ are the extreme eigenvalues of $F$, then
\[
c\big(|y|^2+|\eta|^2\big)
\leq
g_{(x,\xi)}(y,\eta)
\leq
C\big(|y|^2+|\eta|^2\big),
\]
with
\[
c=\min\{\lambda_F^{\min},\lambda_G^{\min}\},
\qquad
C=\max\{\lambda_F^{\max},\lambda_G^{\max}\}.
\]
\end{enumerate}
\end{proposition}

\begin{proof}
Since $G_x$ is positive definite, $T_x$ is invertible and
\[
\|\xi\|_{g_x}=|T_x\xi|
\]
is a norm.  The spectral inclusion gives
\[
\lambda_G^{\min}|\xi|^2
\leq
\xi^\top G_x\xi
\leq
\lambda_G^{\max}|\xi|^2.
\]
Taking square roots proves the Euclidean comparison.

For the comparison between different base points, use the previous estimate twice:
\[
\|\xi\|_{g_x}
\leq
(\lambda_G^{\max})^{1/2}|\xi|
\leq
\left(\frac{\lambda_G^{\max}}{\lambda_G^{\min}}\right)^{1/2}
\|\xi\|_{g_y}.
\]
Interchanging $x$ and $y$ gives the lower bound.

Finally,
\[
g_{(x,\xi)}(y,\eta)
=
y^\top F y+\eta^\top G_x\eta.
\]
Applying the spectral bounds for $F$ and $G_x$ gives
\[
\lambda_F^{\min}|y|^2+\lambda_G^{\min}|\eta|^2
\leq
g_{(x,\xi)}(y,\eta)
\leq
\lambda_F^{\max}|y|^2+\lambda_G^{\max}|\eta|^2,
\]
which is the stated global comparison.
\end{proof}

\begin{remark}[Consequence for the scope of the anisotropy]
The preceding proposition implies
\[
\|\xi\|_{g_x}\asymp |\xi|
\]
uniformly in $x$.  Therefore any symbol estimate expressed only through powers of $\|\xi\|_{g_x}$ is equivalent, up to uniform constants, to the corresponding estimate with powers of $|\xi|$.  The nontrivial feature of the later construction is not a new global symbolic order, but the $x$-dependent localization in the normalized variable $T_x\xi$ and the derivative losses caused by differentiating $T_x$.
\end{remark}

\begin{proposition}[Derivative bounds for $g_x$ and $\|\xi\|_{g_x}$]\label{prop:gx-deriv-bounds-with-xpowers:full}
Under the hypotheses above, the following estimates hold.

\begin{enumerate}
\item The function $g_x(\xi)=\xi^\top G_x\xi$ is smooth in $(x,\xi)$.  For every multi-index $\alpha$ and every multi-index $\beta$ with $|\beta|\leq 2$, there exists $C_{\alpha,\beta}>0$, depending only on $n$, $|\alpha|$, $|\beta|$, and the constants $c_\gamma$ with $|\gamma|\leq|\alpha|$, such that
\begin{equation}\label{eq:bd:dxalpha-dxibeta-g-xpower:full}
|\partial_x^\alpha\partial_\xi^\beta g_x(\xi)|
\leq
C_{\alpha,\beta}
\langle x\rangle^{|\alpha|}
|\xi|^{\max\{2-|\beta|,0\}},
\qquad
x,\xi\in\mathbb{R}^n.
\end{equation}
Moreover, $\partial_\xi^\beta g_x\equiv0$ for $|\beta|\geq3$.

\item For $\xi\neq0$, the function $\|\xi\|_{g_x}=g_x(\xi)^{1/2}$ is smooth in $(x,\xi)$.  For every pair of multi-indices $\alpha,\beta$, there exists $C'_{\alpha,\beta}>0$, depending only on $n$, $|\alpha|$, $|\beta|$, $\lambda_G^{\min}$, $\lambda_G^{\max}$, and the constants $c_\gamma$ with $|\gamma|\leq|\alpha|$, such that
\begin{equation}\label{eq:bd:dxalpha-dxibeta-norm-xpower:full}
|\partial_x^\alpha\partial_\xi^\beta \|\xi\|_{g_x}|
\leq
C'_{\alpha,\beta}
\langle x\rangle^{|\alpha|}
\|\xi\|_{g_x}^{1-|\beta|},
\qquad
x\in\mathbb{R}^n,\quad \xi\neq0.
\end{equation}
\end{enumerate}
\end{proposition}

\begin{proof}
We first prove the estimate for $g_x(\xi)$.  Since $G$ is smooth and $g_x$ is quadratic in $\xi$, smoothness is immediate and all $\xi$-derivatives of order at least $3$ vanish.

If $|\beta|=0$, then
\[
\partial_x^\alpha g_x(\xi)
=
\xi^\top(\partial_x^\alpha G_x)\xi,
\]
so
\[
|\partial_x^\alpha g_x(\xi)|
\leq
\|\partial_x^\alpha G_x\|_{\mathrm{op}}|\xi|^2
\leq
c_\alpha\langle x\rangle^{|\alpha|}|\xi|^2.
\]
If $|\beta|=1$, then $\partial_\xi g_x(\xi)=2G_x\xi$, and therefore
\[
|\partial_x^\alpha\partial_\xi g_x(\xi)|
\leq
2\|\partial_x^\alpha G_x\|_{\mathrm{op}}|\xi|
\leq
2c_\alpha\langle x\rangle^{|\alpha|}|\xi|.
\]
If $|\beta|=2$, then $\partial_\xi^2g_x=2G_x$, hence
\[
|\partial_x^\alpha\partial_\xi^2g_x(\xi)|
\leq
2\|\partial_x^\alpha G_x\|_{\mathrm{op}}
\leq
2c_\alpha\langle x\rangle^{|\alpha|}.
\]
This proves \eqref{eq:bd:dxalpha-dxibeta-g-xpower:full}.

We now prove the estimate for $\|\xi\|_{g_x}$.  Set
\[
\psi(v):=|v|,
\qquad
v(x,\xi):=T_x\xi,
\]
so that $\|\xi\|_{g_x}=\psi(v(x,\xi))$.  On $\mathbb{R}^n\setminus\{0\}$, the function $\psi$ is smooth and homogeneous of degree $1$.  Hence its $r$-th Fréchet derivative is homogeneous of degree $1-r$, and there exists $B_r=B_r(n)>0$ such that
\begin{equation}\label{eq:Drpsi-bound}
\|D^r\psi(v)\|_{\mathrm{op}}
\leq
B_r |v|^{1-r},
\qquad
v\neq0.
\end{equation}

By Lemma~\ref{lem:matrix-sqrt-deriv}, for every multi-index $\gamma$,
\[
\|\partial_x^\gamma T_x\|_{\mathrm{op}}
\leq
C_\gamma\langle x\rangle^{|\gamma|}.
\]
Moreover, uniform ellipticity gives
\[
|\xi|
\leq
(\lambda_G^{\min})^{-1/2}\|\xi\|_{g_x}.
\]

Apply the iterated chain rule to $\psi(T_x\xi)$.  Each term in $\partial_x^\alpha\partial_\xi^\beta\psi(T_x\xi)$ is a finite product obtained by evaluating $D^r\psi(T_x\xi)$ on $r$ factors of the following two types:
\[
(\partial_x^\gamma T_x)\xi,
\qquad
(\partial_x^\gamma T_x)e_j,
\]
where $\gamma$ may be zero in the second case and $e_j$ is a coordinate vector.  The factors of the second type account for the $\xi$-derivatives, and their number is exactly $|\beta|$.  The total number of $x$-derivatives falling on the factors is at most $|\alpha|$.

For a typical term $\mathcal{T}$ in this expansion, suppose that $r$ factors occur and that exactly $|\beta|$ of them are of the second type.  Using \eqref{eq:Drpsi-bound}, Lemma~\ref{lem:matrix-sqrt-deriv}, and the elliptic comparison, we get
\[
\begin{aligned}
|\mathcal{T}|
&\leq
C
\|T_x\xi\|^{1-r}
\langle x\rangle^{|\alpha|}
|\xi|^{r-|\beta|} \\
&\leq
C
\langle x\rangle^{|\alpha|}
\|\xi\|_{g_x}^{1-r}
\|\xi\|_{g_x}^{r-|\beta|} \\
&=
C
\langle x\rangle^{|\alpha|}
\|\xi\|_{g_x}^{1-|\beta|}.
\end{aligned}
\]
There are only finitely many terms, with a number depending on $\alpha$ and $\beta$.  Summing them gives
\[
|\partial_x^\alpha\partial_\xi^\beta \|\xi\|_{g_x}|
\leq
C'_{\alpha,\beta}
\langle x\rangle^{|\alpha|}
\|\xi\|_{g_x}^{1-|\beta|},
\qquad
\xi\neq0.
\]
This proves \eqref{eq:bd:dxalpha-dxibeta-norm-xpower:full}.
\end{proof}

\subsection{Quantization conventions}

We use Weyl quantization with the normalization
\[
\operatorname{Op}^w(a)u(x)
=
(2\pi)^{-n}
\int_{\mathbb{R}^n}
\int_{\mathbb{R}^n}
e^{i(x-y)\cdot\xi}
a\left(\frac{x+y}{2},\xi\right)
u(y)\,\mathrm{d}y\,\mathrm{d}\xi.
\]
Composition of Weyl operators is written in terms of the Moyal product:
\[
\operatorname{Op}^w(a)\operatorname{Op}^w(b)
=
\operatorname{Op}^w(a\# b).
\]
The finite-order expansions and remainders needed for the localized blocks are recalled in Appendix~\ref{sec:appendix}.

%% file: 30_decomp_dyadic.tex
\section{Dyadic microlocal decomposition}
\label{subsec:dyadic-decomposition}

This section constructs the microlocal partition used in the sequel.  The
decomposition remains dyadic in the normalized fiber variable
\[
\zeta=T_x\xi,
\]
in the sense that the centers are organized in the Euclidean annuli
\[
\mathcal C_k
=
\{\zeta\in\mathbb{R}^n:2^k\leq|\zeta|<2^{k+1}\}.
\]
The width of the individual blocks, however, is chosen on the mesoscopic
scale
\[
R_k:=2^{k/2},
\qquad
R_k^2=2^k.
\]
Thus the block diameter is strictly smaller than the dyadic radius,
\[
R_k\ll 2^k
\qquad
\text{as }k\to\infty,
\]
while the reciprocal spatial scale is $R_k^{-1}$.  This choice balances
the variation of $T_x$ across the natural spatial scale of a block and is
the scale used below in the almost-orthogonality estimates.

Since $T_x$ is uniformly elliptic, dyadic size in $\zeta$ is uniformly
equivalent to dyadic size in $\xi$.  The dependence of $T_x$ on $x$
produces nontrivial derivatives of the cutoffs; the mesoscopic
normalization makes these derivatives scale as
\[
\partial_x\sim R_k,
\qquad
\partial_\xi\sim R_k^{-1}.
\]

The construction is made for $k\geq0$.  A fixed low-frequency cutoff is
added separately in order to avoid accumulation of negative dyadic
scales at the zero section.

\subsection{Mesoscopic nets in normalized dyadic annuli}

For $k\in\mathbb{N}_0$ define
\[
\mathcal C_k
:=
\{\zeta\in\mathbb{R}^n:2^k\leq|\zeta|<2^{k+1}\},
\qquad
R_k:=2^{k/2}.
\]
For each $k$, choose a maximal $R_k/2$-separated set
\[
\{\zeta_{j,k}\}_{j\in J_k}\subset\mathcal C_k.
\]
Thus
\[
|\zeta_{j,k}-\zeta_{j',k}|
\geq
\frac{R_k}{2}
\qquad
\text{for }j\neq j',
\]
and maximality implies that for every $\zeta\in\mathcal C_k$ there exists
$j\in J_k$ such that
\[
|\zeta-\zeta_{j,k}|
\leq
\frac{R_k}{2}.
\]
Equivalently,
\[
\mathcal C_k
\subset
\bigcup_{j\in J_k}
\overline B(\zeta_{j,k},R_k/2),
\]
whereas the balls
\[
B(\zeta_{j,k},R_k/4),
\qquad
j\in J_k,
\]
are pairwise disjoint.

The associated fiber centers in the original $\xi$ variable are
\[
\theta_{j,k}(x):=T_x^{-1}\zeta_{j,k}.
\]
Then
\[
T_x\theta_{j,k}(x)=\zeta_{j,k}.
\]
By uniform ellipticity of $T_x$, there exist constants $c,C>0$,
depending only on $\lambda_G^{\min}$ and $\lambda_G^{\max}$, such that
\[
cR_k^2
\leq
|\theta_{j,k}(x)|
\leq
CR_k^2
\]
for every $x$, $j$, and $k$.

The packing estimate for the mesoscopic nets is
\[
\#J_k
\leq
C_n
\left(
\frac{2^k}{R_k}
\right)^n
=
C_nR_k^n
=
C_n2^{kn/2}.
\]
Indeed, the pairwise disjoint balls
$B(\zeta_{j,k},R_k/4)$ are contained in a Euclidean ball of radius
$C2^k$, and the estimate follows by comparison of volumes.

\subsection{Preliminary cutoffs and the low-frequency block}

Fix a nonnegative function $\varphi\in C_c^\infty(\mathbb{R}^n)$ such
that
\[
\varphi(\eta)=1
\quad\text{for }|\eta|\leq\frac12,
\qquad
\operatorname{supp}\varphi
\subset
\{|\eta|\leq1\}.
\]
Also fix a nonnegative function
$\varphi_0\in C_c^\infty(\mathbb{R}^n)$ such that
\[
\varphi_0(\eta)=1
\quad\text{for }|\eta|\leq1,
\qquad
\operatorname{supp}\varphi_0
\subset
\{|\eta|\leq2\}.
\]
For $\ell\in\mathbb{N}_0$ define
\[
M_\ell(\varphi)
:=
\max_{|\alpha|\leq\ell}
\sup_{\eta\in\mathbb{R}^n}
|\partial^\alpha\varphi(\eta)|,
\qquad
M_\ell(\varphi_0)
:=
\max_{|\alpha|\leq\ell}
\sup_{\eta\in\mathbb{R}^n}
|\partial^\alpha\varphi_0(\eta)|.
\]

For each pair $(j,k)$ set
\[
u_{j,k}(x,\xi)
:=
\frac{T_x\xi-\zeta_{j,k}}{R_k}
=
\frac{T_x(\xi-\theta_{j,k}(x))}{R_k},
\]
and define
\[
\chi_{j,k}(x,\xi)
:=
\varphi\bigl(u_{j,k}(x,\xi)\bigr).
\]
Thus
\[
\operatorname{supp}\chi_{j,k}(x,\cdot)
\subset
\left\{
\xi:
|T_x\xi-\zeta_{j,k}|
\leq
R_k
\right\},
\]
and
\[
\chi_{j,k}(x,\xi)=1
\qquad
\text{whenever }
|T_x\xi-\zeta_{j,k}|
\leq
\frac{R_k}{2}.
\]

On the support of $\chi_{j,k}$ one has
\[
\langle T_x\xi\rangle
\asymp
R_k^2
\]
uniformly in $j,k$.  Hence, by uniform ellipticity,
\begin{equation}
\langle\xi\rangle
\asymp
R_k^2
\qquad
\text{on }\operatorname{supp}\chi_{j,k},
\end{equation}
with constants independent of $j$ and $k$.

The low-frequency cutoff remains fixed:
\[
\chi_0(x,\xi)
:=
\varphi_0(T_x\xi).
\]
Hence
\[
\chi_0(x,\xi)=1
\quad\text{if }|T_x\xi|\leq1,
\qquad
\operatorname{supp}\chi_0
\subset
\{(x,\xi):|T_x\xi|\leq2\}.
\]

\begin{lemma}[Uniform overlap]\label{lem:besicovitch-coronas-radial-full}
There exists a constant $L=L(n)$ such that
\[
\sup_{x,\xi}
\left(
\mathbf{1}_{\operatorname{supp}\chi_0}(x,\xi)
+
\sum_{k\geq0}\sum_{j\in J_k}
\mathbf{1}_{\operatorname{supp}\chi_{j,k}(x,\cdot)}(\xi)
\right)
\leq
L.
\]
The constant depends only on $n$ and on the fixed support radii of
$\varphi$ and $\varphi_0$.  Moreover,
\[
\#J_k\leq C_nR_k^n
\]
for every $k\geq0$.
\end{lemma}

\begin{proof}
Fix $x\in\mathbb{R}^n$ and $\xi\in\mathbb{R}^n$, and write
\[
\zeta=T_x\xi.
\]
The low-frequency cutoff contributes at most one term.

Suppose
\[
\xi\in\operatorname{supp}\chi_{j,k}(x,\cdot).
\]
Then
\[
|\zeta-\zeta_{j,k}|
\leq
R_k.
\]
Since
\[
\zeta_{j,k}\in\mathcal C_k,
\]
we have
\[
2^k-R_k
\leq
|\zeta|
<
2^{k+1}+R_k.
\]
For $k\geq2$,
\[
R_k=2^{k/2}\leq2^{k-1},
\]
and therefore
\[
2^{k-1}
\leq
|\zeta|
<
2^{k+2}.
\]
Hence only uniformly many values of $k\geq2$ can occur for a fixed
$\zeta$.  The remaining values $k=0,1$ form a fixed finite set.
Consequently, the number of dyadic scales contributing at any point is
uniformly bounded.

For fixed $k$, let
\[
S_{x,k}(\xi)
:=
\{
j\in J_k:
\xi\in\operatorname{supp}\chi_{j,k}(x,\cdot)
\}.
\]
If $j\in S_{x,k}(\xi)$, then
\[
\zeta_{j,k}\in B(\zeta,R_k).
\]
The balls
\[
B(\zeta_{j,k},R_k/4),
\qquad
j\in S_{x,k}(\xi),
\]
are pairwise disjoint and are contained in
\[
B(\zeta,5R_k/4).
\]
Hence
\[
\#S_{x,k}(\xi)
\operatorname{Vol}(B(0,R_k/4))
\leq
\operatorname{Vol}(B(0,5R_k/4)),
\]
so that
\[
\#S_{x,k}(\xi)\leq5^n.
\]
Combining this estimate with the uniform bound on the number of
admissible scales proves the overlap assertion.

Finally, the balls $B(\zeta_{j,k},R_k/4)$, $j\in J_k$, are pairwise
disjoint and contained in
\[
B(0,2^{k+1}+R_k/4)
\subset
B(0,C2^k).
\]
Comparison of volumes gives
\[
\#J_k
\leq
C_n
\left(
\frac{2^k}{R_k}
\right)^n
=
C_nR_k^n.
\]
\end{proof}

\subsection{The normalizing sum and the partition}

Define
\[
\Sigma(x,\xi)
:=
\chi_0(x,\xi)
+
\sum_{k\geq0}\sum_{j\in J_k}
\chi_{j,k}(x,\xi).
\]
The sum is pointwise finite by
Lemma~\ref{lem:besicovitch-coronas-radial-full}.

Moreover, $\Sigma$ is uniformly bounded below.  Indeed, if
\[
|T_x\xi|\leq1,
\]
then
\[
\chi_0(x,\xi)=1.
\]
If
\[
|T_x\xi|\geq1,
\]
then, by the half-open convention in the definition of $\mathcal C_k$,
there is a unique $k\in\mathbb{N}_0$ such that
\[
T_x\xi\in\mathcal C_k.
\]
By maximality of the $R_k/2$-separated net, there exists
$j\in J_k$ such that
\[
|T_x\xi-\zeta_{j,k}|
\leq
\frac{R_k}{2}.
\]
For this pair,
\[
\chi_{j,k}(x,\xi)=1.
\]
Thus
\begin{equation}\label{eq:Sigma-lower-bound}
\Sigma(x,\xi)\geq1
\qquad
\text{for every }(x,\xi)\in\mathbb{R}^{2n}.
\end{equation}

Define
\[
\Lambda_0(x,\xi)
:=
\frac{\chi_0(x,\xi)}{\Sigma(x,\xi)}
\]
and, for $k\geq0$,
\[
\Lambda_{j,k}(x,\xi)
:=
\frac{\chi_{j,k}(x,\xi)}{\Sigma(x,\xi)}.
\]
Then
\[
0\leq\Lambda_0\leq1,
\qquad
0\leq\Lambda_{j,k}\leq1,
\]
the family is locally finite, and
\begin{equation}\label{eq:dyadic-partition-unity}
\Lambda_0(x,\xi)
+
\sum_{k\geq0}\sum_{j\in J_k}
\Lambda_{j,k}(x,\xi)
=
1
\qquad
\text{for every }(x,\xi)\in\mathbb{R}^{2n}.
\end{equation}

Furthermore,
\[
\operatorname{supp}\Lambda_{j,k}(x,\cdot)
\subset
\left\{
\xi:
|T_x\xi-\zeta_{j,k}|
\leq
R_k
\right\},
\]
and consequently
\[
\langle\xi\rangle
\asymp
R_k^2
\qquad
\text{on }\operatorname{supp}\Lambda_{j,k},
\]
uniformly in $j$ and $k$.

\begin{remark}
The dyadic radius of the $k$-th annulus is $R_k^2=2^k$, whereas the
individual normalized blocks have radius $R_k=2^{k/2}$.  Thus the
decomposition is dyadic at the level of the annuli and mesoscopic at the
level of the individual blocks.  The fixed block $\Lambda_0$ contains
the low-frequency region and is treated separately by the standard
low-frequency calculus.
\end{remark}

\subsection{Scale-adapted derivative estimates}

We first estimate the unnormalized high-frequency cutoffs.  Since
\[
u_{j,k}(x,\xi)
=
R_k^{-1}(T_x\xi-\zeta_{j,k}),
\]
and $\zeta_{j,k}$ is independent of $x$, one has
\[
\partial_x^\alpha\partial_\xi^\beta
u_{j,k}(x,\xi)
=
\begin{cases}
R_k^{-1}(\partial_x^\alpha T_x)\xi,
&
\beta=0,
\\[4pt]
R_k^{-1}(\partial_x^\alpha T_x)e_\ell,
&
\beta=e_\ell,
\\[4pt]
0,
&
|\beta|\geq2.
\end{cases}
\]
By Lemma~\ref{lem:matrix-sqrt-deriv},
\[
\|\partial_x^\alpha T_x\|_{\mathrm{op}}
\leq
C_\alpha\langle x\rangle^{|\alpha|}.
\]
Therefore
\begin{equation}\label{eq:u-deriv-bound}
|\partial_x^\alpha\partial_\xi^\beta
u_{j,k}(x,\xi)|
\leq
C_{\alpha,\beta}
\langle x\rangle^{|\alpha|}
R_k^{-1}
\langle\xi\rangle^{\mathbf{1}_{\{|\beta|=0\}}}
\end{equation}
whenever $|\beta|\leq1$, while the derivative vanishes for
$|\beta|\geq2$.

On the support of $\chi_{j,k}$,
\[
\langle\xi\rangle\asymp R_k^2,
\]
so \eqref{eq:u-deriv-bound} yields the scale-adapted estimates
\[
|\partial_x^\alpha u_{j,k}(x,\xi)|
\leq
C_\alpha
\langle x\rangle^{|\alpha|}
R_k,
\]
and
\[
|\partial_x^\alpha\partial_{\xi_\ell}
u_{j,k}(x,\xi)|
\leq
C_\alpha
\langle x\rangle^{|\alpha|}
R_k^{-1}
\]
there.

\begin{lemma}[Single-cutoff derivative bound]\label{lem:single-chi-deriv}
For every pair of multi-indices $\alpha,\beta$, there exists
$C_{\alpha,\beta}>0$ such that
\[
|\partial_x^\alpha\partial_\xi^\beta
\chi_{j,k}(x,\xi)|
\leq
C_{\alpha,\beta}
\langle x\rangle^{|\alpha|}
R_k^{|\alpha|-|\beta|}
\]
for all $x,\xi,j,k$ with $k\geq0$.

For the fixed low-frequency cutoff one has
\[
|\partial_x^\alpha\partial_\xi^\beta
\chi_0(x,\xi)|
\leq
C_{\alpha,\beta}
\langle x\rangle^{|\alpha|}.
\]

The constants depend only on finitely many seminorms of $\varphi$ and
$\varphi_0$ and on the constants in
Lemma~\ref{lem:matrix-sqrt-deriv} up to order $|\alpha|$.
\end{lemma}

\begin{proof}
We apply the iterated chain rule to
\[
\chi_{j,k}=\varphi\circ u_{j,k}.
\]
Every nonzero chain-rule term is a derivative of $\varphi$ evaluated at
$u_{j,k}$ multiplied by factors of the form
\[
\partial_x^{\alpha_q}
\partial_\xi^{\beta_q}
u_{j,k},
\qquad
|\beta_q|\leq1.
\]

A factor with $\beta_q=0$ is nontrivial only if it consumes at least one
$x$-derivative, except for the undifferentiated expression, and on the
support of the corresponding derivative it is bounded by
\[
C
\langle x\rangle^{|\alpha_q|}
R_k.
\]
A factor with $|\beta_q|=1$ is bounded by
\[
C
\langle x\rangle^{|\alpha_q|}
R_k^{-1}.
\]
The total number of factors of the first type is at most $|\alpha|$,
whereas the total number of factors carrying one $\xi$-derivative is
exactly $|\beta|$.  Consequently every term is bounded by
\[
C_{\alpha,\beta}
\langle x\rangle^{|\alpha|}
R_k^{|\alpha|-|\beta|}.
\]
The derivatives of $\varphi$ are uniformly bounded, which proves the
high-frequency estimate.

For $\chi_0=\varphi_0(T_x\xi)$, the support condition
\[
|T_x\xi|\leq2
\]
and uniform ellipticity imply $|\xi|\leq C$.  The same chain-rule
argument therefore gives
\[
|\partial_x^\alpha\partial_\xi^\beta
\chi_0(x,\xi)|
\leq
C_{\alpha,\beta}
\langle x\rangle^{|\alpha|}.
\]
\end{proof}

\begin{proposition}[Derivative bounds for $\Sigma$]
\label{prop:Sigma-deriv-bounds:proof}
The function $\Sigma$ is smooth on $\mathbb{R}^{2n}$.  For every pair of
multi-indices $\alpha,\beta$, there exists $C_{\alpha,\beta}>0$ such
that
\[
|\partial_x^\alpha\partial_\xi^\beta
\Sigma(x,\xi)|
\leq
C_{\alpha,\beta}
\langle x\rangle^{|\alpha|}
\langle\xi\rangle^{\frac{|\alpha|-|\beta|}{2}}
\]
for all $(x,\xi)\in\mathbb{R}^{2n}$.

The constant depends only on $n$, finitely many seminorms of $\varphi$
and $\varphi_0$, the structural constants of $G$, and the overlap
constant in Lemma~\ref{lem:besicovitch-coronas-radial-full}.
\end{proposition}

\begin{proof}
Local finiteness follows from
Lemma~\ref{lem:besicovitch-coronas-radial-full}, so $\Sigma$ is smooth.

At every point $(x,\xi)$ only uniformly many terms in the definition of
$\Sigma$ are nonzero.  The same is true after differentiation, since
the derivatives of each cutoff are supported in its original support.

If a high-frequency block of scale $k$ contributes at $(x,\xi)$, then
\[
\langle\xi\rangle\asymp R_k^2.
\]
By Lemma~\ref{lem:single-chi-deriv},
\[
|\partial_x^\alpha\partial_\xi^\beta
\chi_{j,k}(x,\xi)|
\leq
C_{\alpha,\beta}
\langle x\rangle^{|\alpha|}
R_k^{|\alpha|-|\beta|}
\leq
C_{\alpha,\beta}
\langle x\rangle^{|\alpha|}
\langle\xi\rangle^{\frac{|\alpha|-|\beta|}{2}}.
\]
All high-frequency scales contributing at the same point are comparable,
and their number is uniformly bounded.

The low-frequency term is supported where $\langle\xi\rangle\asymp1$,
and therefore satisfies the same estimate after enlarging the
constant.  Summing the uniformly bounded number of nonzero terms proves
the claim.
\end{proof}

We now estimate the reciprocal normalizing factor.

\begin{lemma}[Derivative bounds for the reciprocal]
\label{lem:Sigma-reciprocal-deriv}
For every pair of multi-indices $\alpha,\beta$, there exists
$C_{\alpha,\beta}>0$ such that
\[
|\partial_x^\alpha\partial_\xi^\beta
(\Sigma^{-1})(x,\xi)|
\leq
C_{\alpha,\beta}
\langle x\rangle^{|\alpha|}
\langle\xi\rangle^{\frac{|\alpha|-|\beta|}{2}}
\]
for all $(x,\xi)\in\mathbb{R}^{2n}$.
\end{lemma}

\begin{proof}
The case $\alpha=\beta=0$ follows from
\eqref{eq:Sigma-lower-bound}.  For higher derivatives, apply the
iterated chain rule to
\[
F\circ\Sigma,
\qquad
F(t)=t^{-1}.
\]
Each term has the form
\[
F^{(r)}(\Sigma)
\prod_{\nu=1}^r
\partial_x^{\alpha_\nu}
\partial_\xi^{\beta_\nu}\Sigma,
\]
where
\[
\sum_{\nu=1}^r\alpha_\nu=\alpha,
\qquad
\sum_{\nu=1}^r\beta_\nu=\beta.
\]
Since $\Sigma\geq1$,
\[
|F^{(r)}(\Sigma)|\leq r!.
\]
Using Proposition~\ref{prop:Sigma-deriv-bounds:proof},
the product is bounded by
\[
C_{\alpha,\beta}
\langle x\rangle^{\sum_\nu|\alpha_\nu|}
\langle\xi\rangle^{
\frac12
\sum_\nu(|\alpha_\nu|-|\beta_\nu|)
}
=
C_{\alpha,\beta}
\langle x\rangle^{|\alpha|}
\langle\xi\rangle^{\frac{|\alpha|-|\beta|}{2}}.
\]
Summing the finitely many chain-rule terms proves the assertion.
\end{proof}

\begin{proposition}[Derivative bounds for the normalized cutoffs]
\label{prop:Lambda-deriv-bounds:full}
Let
\[
\Lambda_0=\chi_0\Sigma^{-1},
\qquad
\Lambda_{j,k}
=
\chi_{j,k}\Sigma^{-1}
\quad
(k\geq0).
\]
For every pair of multi-indices $\alpha,\beta$, there exists
$C_{\alpha,\beta}>0$ such that
\begin{equation}\label{eq:bd:Lambda-general-xpower:full}
|\partial_x^\alpha\partial_\xi^\beta
\Lambda_{j,k}(x,\xi)|
\leq
C_{\alpha,\beta}
\langle x\rangle^{|\alpha|}
R_k^{|\alpha|-|\beta|}
\end{equation}
for all $(x,\xi)\in\mathbb{R}^{2n}$ and all $j,k$ with $k\geq0$.

Equivalently, one has the global estimate
\[
|\partial_x^\alpha\partial_\xi^\beta
\Lambda_{j,k}(x,\xi)|
\leq
C_{\alpha,\beta}
\langle x\rangle^{|\alpha|}
\langle\xi\rangle^{\frac{|\alpha|-|\beta|}{2}}
\mathbf{1}_{\operatorname{supp}\chi_{j,k}}(x,\xi).
\]

For the low-frequency block,
\[
|\partial_x^\alpha\partial_\xi^\beta
\Lambda_0(x,\xi)|
\leq
C_{\alpha,\beta}
\langle x\rangle^{|\alpha|}.
\]
All constants are independent of $j$ and $k$.
\end{proposition}

\begin{proof}
By Leibniz' rule,
\[
\partial_x^\alpha\partial_\xi^\beta
\Lambda_{j,k}
=
\sum_{\substack{\alpha_1+\alpha_2=\alpha\\
\beta_1+\beta_2=\beta}}
\binom{\alpha}{\alpha_1}
\binom{\beta}{\beta_1}
\bigl(
\partial_x^{\alpha_1}\partial_\xi^{\beta_1}
\chi_{j,k}
\bigr)
\bigl(
\partial_x^{\alpha_2}\partial_\xi^{\beta_2}
\Sigma^{-1}
\bigr).
\]
Every summand is supported in
$\operatorname{supp}\chi_{j,k}$.  On this set,
\[
\langle\xi\rangle\asymp R_k^2.
\]
Using Lemma~\ref{lem:single-chi-deriv} and
Lemma~\ref{lem:Sigma-reciprocal-deriv}, we obtain
\[
\begin{aligned}
&
\left|
\bigl(
\partial_x^{\alpha_1}\partial_\xi^{\beta_1}
\chi_{j,k}
\bigr)
\bigl(
\partial_x^{\alpha_2}\partial_\xi^{\beta_2}
\Sigma^{-1}
\bigr)
\right|
\\
&\qquad
\leq
C
\langle x\rangle^{|\alpha_1|+|\alpha_2|}
R_k^{|\alpha_1|-|\beta_1|}
R_k^{|\alpha_2|-|\beta_2|}
\\
&\qquad
=
C
\langle x\rangle^{|\alpha|}
R_k^{|\alpha|-|\beta|}.
\end{aligned}
\]
Summing the finitely many Leibniz terms proves
\eqref{eq:bd:Lambda-general-xpower:full}.

The equivalent estimate in terms of $\langle\xi\rangle$ follows again
from
\[
\langle\xi\rangle\asymp R_k^2
\qquad
\text{on }\operatorname{supp}\chi_{j,k}.
\]
The proof for $\Lambda_0$ is identical, using the fixed-scale bounds for
$\chi_0$ and the fact that its support lies in a bounded normalized
frequency region.
\end{proof}

\begin{remark}[Mesoscopic derivative balance]
The normalized high-frequency blocks satisfy the scale-adapted estimate
\[
|\partial_x^\alpha\partial_\xi^\beta
\Lambda_{j,k}(x,\xi)|
\leq
C_{\alpha,\beta}
\langle x\rangle^{|\alpha|}
R_k^{|\alpha|-|\beta|}.
\]
Thus one $x$-derivative costs one factor $R_k$, whereas one
$\xi$-derivative gains one factor $R_k^{-1}$.  Since
\[
R_k^2=2^k
\]
is the dyadic frequency scale, the corresponding local phase-space
metric is
\[
R_k^2|\mathrm{d}x|^2
+
R_k^{-2}|\mathrm{d}\xi|^2.
\]
On a fixed spatial compact, the factors
$\langle x\rangle^{|\alpha|}$ are absorbed into the constants, and the
high-frequency cutoffs satisfy uniformly
\[
|\partial_x^\alpha\partial_\xi^\beta
\Lambda_{j,k}(x,\xi)|
\leq
C_{\alpha,\beta,K}
R_k^{|\alpha|-|\beta|}.
\]
This is the scale used in the subsequent kernel localization and
almost-orthogonality estimates.
\end{remark}

%% file: 40_quanti_dyadic.tex
\section{Quantization and mesoscopic recombination}
\label{sec:quantization-local-global}

This section contains the quantitative core of the paper.  The mesoscopic
scale
\[
R_k:=2^{k/2},
\qquad
R_k^2=2^k,
\]
introduced in Section~\ref{subsec:dyadic-decomposition} is used throughout.
The first result records both a conservative finite-seminorm estimate in the
global symbol class and the sharper scale-adapted estimate on each
high-frequency block.  The latter is the estimate relevant for Weyl
quantization and almost orthogonality.  In particular, the artificial band
loss produced by the former fixed-radius decomposition no longer appears in
the Sobolev-conjugated block estimates.

Throughout the section, $\operatorname{Op}^w$ denotes Weyl quantization with
the convention fixed in Section~\ref{sec:metricas}.  If
$K\subset\mathbb{R}^n$ is compact, we use the localized seminorm
\[
\mathcal S_{\imath,\ell}(a;K)
:=
\sup_{\substack{|\alpha|\leq\imath\\|\beta|\leq\ell}}
\sup_{(x,\xi)\in K\times\mathbb{R}^n}
\frac{|\partial_x^\alpha\partial_\xi^\beta a(x,\xi)|}
{\langle x\rangle^{m_1-|\alpha|}
 \langle\xi\rangle^{m_2-|\beta|}}.
\]
The corresponding global seminorm is denoted by
$\mathcal S_{\imath,\ell}(a)$.

\subsection{Finite-seminorm stability under mesoscopic localization}

\begin{proposition}[Finite-seminorm localization at the mesoscopic scale]
\label{prop:mult-by-Lambda-preserves-S-anisotropic:full}
Let $a\in S^{m_1,m_2}$ and fix $\imath,\ell\in\mathbb{N}_0$.  Set
\[
N_x:=2\imath,
\qquad
N_\xi:=\left\lceil\frac{\imath+\ell}{2}\right\rceil.
\]
Then there exists $C_{\imath,\ell}>0$, independent of $(j,k)$, such that
\[
\|a\Lambda_{j,k}\|_{\imath,\ell}^{(m_1+N_x,m_2+N_\xi)}
\leq
C_{\imath,\ell}
\|a\|_{\imath,\ell}^{(m_1,m_2)}.
\]
Equivalently, for $|\alpha|\leq\imath$ and $|\beta|\leq\ell$,
\[
|\partial_x^\alpha\partial_\xi^\beta(a\Lambda_{j,k})(x,\xi)|
\leq
C_{\imath,\ell}
\|a\|_{\imath,\ell}^{(m_1,m_2)}
\langle x\rangle^{m_1+N_x-|\alpha|}
\langle\xi\rangle^{m_2+N_\xi-|\beta|}.
\]
The same conservative estimate holds for the fixed low-frequency block
$a\Lambda_0$.

Moreover, if $K\subset\mathbb{R}^n$ is compact, then the high-frequency
blocks satisfy the sharper estimate
\begin{equation}\label{eq:localized-symbol-mesoscopic}
|\partial_x^\alpha\partial_\xi^\beta
(a\Lambda_{j,k})(x,\xi)|
\leq
C_{\alpha,\beta,K}
\mathcal S_{\imath,\ell}(a;K)
R_k^{2m_2+|\alpha|-|\beta|}
\end{equation}
for $x\in K$, $|\alpha|\leq\imath$, $|\beta|\leq\ell$, and
$(x,\xi)\in\operatorname{supp}\Lambda_{j,k}$.  The constant in
\eqref{eq:localized-symbol-mesoscopic} is independent of $(j,k)$.
\end{proposition}

\begin{proof}
Fix $|\alpha|\leq\imath$ and $|\beta|\leq\ell$.  By Leibniz' rule,
\[
\partial_x^\alpha\partial_\xi^\beta(a\Lambda_{j,k})
=
\sum_{\substack{\alpha_1+\alpha_2=\alpha\\
\beta_1+\beta_2=\beta}}
\binom{\alpha}{\alpha_1}
\binom{\beta}{\beta_1}
(\partial_x^{\alpha_1}\partial_\xi^{\beta_1}a)
(\partial_x^{\alpha_2}\partial_\xi^{\beta_2}\Lambda_{j,k}).
\]
Proposition~\ref{prop:Lambda-deriv-bounds:full} gives
\[
|\partial_x^{\alpha_2}\partial_\xi^{\beta_2}
\Lambda_{j,k}(x,\xi)|
\leq
C_{\alpha_2,\beta_2}
\langle x\rangle^{|\alpha_2|}
\langle\xi\rangle^{\frac{|\alpha_2|-|\beta_2|}{2}}
\mathbf 1_{\operatorname{supp}\chi_{j,k}}(x,\xi).
\]
Combining this with the symbol estimate for $a$, each Leibniz term is
bounded by
\[
C
\|a\|_{\imath,\ell}^{(m_1,m_2)}
\langle x\rangle^{m_1-|\alpha_1|+|\alpha_2|}
\langle\xi\rangle^{m_2-|\beta_1|
+\frac{|\alpha_2|-|\beta_2|}{2}}.
\]
After division by
\[
\langle x\rangle^{m_1+N_x-|\alpha|}
\langle\xi\rangle^{m_2+N_\xi-|\beta|},
\]
the residual powers are
\[
\langle x\rangle^{2|\alpha_2|-N_x}
\langle\xi\rangle^{\frac{|\alpha_2|+|\beta_2|}{2}-N_\xi}.
\]
The choices
\[
N_x=2\imath,
\qquad
N_\xi=\left\lceil\frac{\imath+\ell}{2}\right\rceil
\]
make both exponents nonpositive.  This proves the conservative global
finite-seminorm estimate.  The argument for $\Lambda_0$ is identical,
using the fixed low-frequency bounds in
Proposition~\ref{prop:Lambda-deriv-bounds:full}.

For the scale-adapted estimate, restrict $x$ to $K$.  On
$\operatorname{supp}\Lambda_{j,k}$,
\[
\langle\xi\rangle\asymp R_k^2.
\]
Hence
\[
|\partial_x^{\alpha_1}\partial_\xi^{\beta_1}a(x,\xi)|
\leq
C_K\mathcal S_{\imath,\ell}(a;K)
R_k^{2m_2-2|\beta_1|},
\]
while
\[
|\partial_x^{\alpha_2}\partial_\xi^{\beta_2}
\Lambda_{j,k}(x,\xi)|
\leq
C_{\alpha_2,\beta_2,K}
R_k^{|\alpha_2|-|\beta_2|}.
\]
Therefore every Leibniz term is bounded by
\[
C_K\mathcal S_{\imath,\ell}(a;K)
R_k^{2m_2+|\alpha_2|-2|\beta_1|-|\beta_2|}
\leq
C_K\mathcal S_{\imath,\ell}(a;K)
R_k^{2m_2+|\alpha|-|\beta|},
\]
which proves \eqref{eq:localized-symbol-mesoscopic}.
\end{proof}

\begin{remark}[Global finite loss versus mesoscopic block control]
\label{rem:no-global-symbol-class-from-finite-loss}
The first part of
Proposition~\ref{prop:mult-by-Lambda-preserves-S-anisotropic:full} is a
conservative finite-seminorm statement in the global class
$S^{m_1,m_2}$.  The second part is sharper on each high-frequency block:
one $x$-derivative costs one factor $R_k$, whereas one $\xi$-derivative
gains one factor $R_k^{-1}$.  The local Weyl and almost-orthogonality
arguments below use this scale-adapted estimate rather than converting the
cutoff derivatives into an artificial positive symbolic order.
\end{remark}

\subsection{Local Weyl bounds}

\begin{lemma}[High-frequency support of the mesoscopic blocks]
\label{lem:dyadic-support-xi}
There exist $k_\ast\in\mathbb{N}_0$ and constants $0<c<C<\infty$,
depending only on the ellipticity constants of $G_x$, such that for every
$k\geq k_\ast$,
\[
\operatorname{supp}\Lambda_{j,k}
\subset
\{(x,\xi):cR_k^2\leq|\xi|\leq CR_k^2\}.
\]
Equivalently,
\[
\langle\xi\rangle\asymp R_k^2=2^k
\]
on $\operatorname{supp}\Lambda_{j,k}$, uniformly in $x,j,k$.
\end{lemma}

\begin{proof}
On $\operatorname{supp}\Lambda_{j,k}$ one has
\[
|T_x\xi-\zeta_{j,k}|\leq R_k,
\qquad
R_k^2\leq|\zeta_{j,k}|<2R_k^2.
\]
Thus
\[
R_k^2-R_k
\leq
|T_x\xi|
\leq
2R_k^2+R_k.
\]
Choose $k_\ast$ so that $R_k\geq2$ for $k\geq k_\ast$.  Then
\[
\frac12R_k^2
\leq
|T_x\xi|
\leq
\frac52R_k^2.
\]
Uniform ellipticity of $T_x$ gives the claimed comparison with $|\xi|$.
\end{proof}

\begin{lemma}[Local Sobolev reduction to Calder\'on--Vaillancourt]
\label{lem:local-sobolev-cv-reduction}
Let $K\subset\mathbb{R}^n$ be compact, let
$\chi,\chi'\in C_c^\infty(\mathbb{R}^n)$ with
$\operatorname{supp}\chi,\operatorname{supp}\chi'\subset K$, and let
$s,m\in\mathbb{R}$.  There exist integers $A=A(n,s,m)$ and
$B=B(n,s,m)$ such that, for every symbol $p$ satisfying
\[
M_{A,B}^{(m)}(p;K)
:=
\sup_{\substack{|\alpha|\leq A\\|\beta|\leq B}}
\sup_{(x,\xi)\in K\times\mathbb{R}^n}
\langle\xi\rangle^{-m+|\beta|}
|\partial_x^\alpha\partial_\xi^\beta p(x,\xi)|
<\infty,
\]
one has
\[
\|\chi\operatorname{Op}^w(p)\chi'\|_{H^s\to H^{s-m}}
\leq
C M_{A,B}^{(m)}(p;K_1),
\]
where $K_1$ is a fixed compact neighborhood of $K$ depending only on the
cutoffs.
\end{lemma}

\begin{proof}
This is the local Sobolev form of the Calder\'on--Vaillancourt estimate.
After conjugation by $\langle D\rangle^{s-m}$ and
$\langle D\rangle^{-s}$, a finite Weyl expansion reduces the operator to
an order-zero symbol whose finite seminorms are controlled by
$M_{A,B}^{(m)}(p;K_1)$.  Theorem~\ref{thm:CV-local} then gives the stated
bound.
\end{proof}

\begin{corollary}[Scale-adapted derivative estimate on one block]
\label{cor:weyl-operator-bounds-with-loss:local-theta}
Let $a\in S^{m_1,m_2}$, let $K\subset\mathbb{R}^n$ be compact, and let
$\chi,\chi'\in C_c^\infty(K)$.  For every finite $A,B$ there exists a
constant $C_{A,B,K}$ such that, for $k\geq k_\ast$,
\[
|\partial_x^\alpha\partial_\xi^\beta
(a\Lambda_{j,k})(x,\xi)|
\leq
C_{A,B,K}\mathcal S_{A,B}(a;K)
R_k^{2m_2+|\alpha|-|\beta|}
\]
for $|\alpha|\leq A$, $|\beta|\leq B$, uniformly in $j,k$.
\end{corollary}

\begin{proof}
This is the localized estimate
\eqref{eq:localized-symbol-mesoscopic}.
\end{proof}

\begin{definition}[Sobolev-conjugated blocks]
\label{def:sobolev-conjugated-blocks}
Let $\psi,\psi'\in C_c^\infty(\mathbb{R}^n)$ be fixed spatial cutoffs and
set, for $k\geq0$,
\[
T_{j,k}
:=
\psi\operatorname{Op}^w(a\Lambda_{j,k})\psi'
:
H^s\to H^{s-m_2}.
\]
The associated $L^2$ block is
\[
\widetilde T_{j,k}
:=
\langle D\rangle^{s-m_2}
T_{j,k}
\langle D\rangle^{-s}
:
L^2\to L^2.
\]
For $\lambda=(j,k)$ we write
\[
R_\lambda:=R_k,
\qquad
\zeta_\lambda:=\zeta_{j,k}.
\]
For two high-frequency indices $\lambda,\mu$ define
\begin{equation}\label{eq:mesoscopic-block-distance}
\operatorname{dist}_{\mathrm{mes}}(\lambda,\mu)
:=
\frac{|\zeta_\lambda-\zeta_\mu|}
{R_\lambda+R_\mu}.
\end{equation}
\end{definition}

\begin{lemma}[Sobolev conjugation and mesoscopic kernel localization]
\label{lem:mesoscopic-kernel-localization}
Let $a\in S^{m_1,m_2}$ and fix $s\in\mathbb{R}$ and
$\psi,\psi'\in C_c^\infty(\mathbb{R}^n)$.  Let
$\widetilde T_\lambda$ be as in
Definition~\ref{def:sobolev-conjugated-blocks}.  For every finite $A,B,L\in\mathbb{N}_0$, there exists a fixed compact set
$K_1\subset\mathbb{R}^n$, independent of $\lambda$, such that the exact
Weyl symbol $b_\lambda$ of $\widetilde T_\lambda$ satisfies
\begin{equation}\label{eq:mesoscopic-conjugated-symbol}
|\partial_x^\alpha\partial_\xi^\beta b_\lambda(x,\xi)|
\leq
C_{A,B,L}
\mathcal S_{D,D}(a;K_1)
R_\lambda^{|\alpha|-|\beta|}
\left(
1+\operatorname{dist}_{g_\lambda}
\big((x,\xi),\Omega_\lambda\big)
\right)^{-L}
\end{equation}
for $|\alpha|\leq A$ and $|\beta|\leq B$, where
\[
\Omega_\lambda
:=
\left\{
(y,\eta):
y\in K_1,
\ |T_y\eta-\zeta_\lambda|\leq C_1R_\lambda
\right\}
\]
and
\[
\operatorname{dist}_{g_\lambda}
\big((x,\xi),\Omega_\lambda\big)
:=
\inf_{(y,\eta)\in\Omega_\lambda}
\left(
R_\lambda|x-y|
+
R_\lambda^{-1}|\xi-\eta|
\right).
\]
The constants are independent of $j$ and $k$.

Let $S_x:=T_x^{-1}$.  On every fixed compact neighborhood
$K_2\supset K_1$, the Schwartz kernel of $\widetilde T_\lambda$ can be
written exactly as
\begin{equation}\label{eq:mesoscopic-kernel-form}
K_\lambda(x,y)
=
R_\lambda^n
\exp\!\left(
 i(x-y)\cdot
 S_{(x+y)/2}\zeta_\lambda
\right)
A_\lambda\!\left(
\frac{x+y}{2},
R_\lambda(x-y)
\right),
\end{equation}
where $D=D(n,s,m_2,A,B,L)$ is finite.  For the Cotlar--Stein
application with the fixed exponent
\[
M_{\operatorname{CS}}=2n+2,
\]
one may take
\[
D=D_{\operatorname{CS}},
\]
with $D_{\operatorname{CS}}$ given explicitly in
\eqref{eq:CS-derivative-budget}; see
Proposition~\ref{prop:mesoscopic-cotlar-derivative-budget}.

\begin{equation}\label{eq:mesoscopic-amplitude-decay}
|\partial_m^\alpha\partial_u^\beta A_\lambda(m,u)|
\leq
C_{\alpha,\beta,N}
\mathcal S_{D,D}(a;K_1)
R_\lambda^{|\alpha|}
\langle u\rangle^{-N}
\end{equation}
for $m\in K_2$ and every prescribed finite $\alpha,\beta,N$.  The part of
the kernel for which $(x+y)/2\notin K_2$ satisfies estimates of arbitrary
order obtained from \eqref{eq:mesoscopic-conjugated-symbol} and can be
absorbed into the same bounds after increasing $N$.
\end{lemma}

\begin{proof}
Write
\[
p_\lambda:=a\Lambda_\lambda,
\qquad
w_+(\xi):=\langle\xi\rangle^{s-m_2},
\qquad
w_-(\xi):=\langle\xi\rangle^{-s}.
\]
The exact Weyl symbol of the conjugated block is
\[
b_\lambda
=
w_+\#\psi\#p_\lambda\#\psi'\#w_-.
\]
Choose $K_1$ containing the spatial supports of the cutoffs and all
midpoints generated by them.  Proposition~\ref{prop:mult-by-Lambda-preserves-S-anisotropic:full}
gives, on $K_1$,
\[
|\partial_x^\alpha\partial_\xi^\beta p_\lambda(x,\xi)|
\leq
C
\mathcal S_{D,D}(a;K_1)
R_\lambda^{2m_2+|\alpha|-|\beta|}.
\]
Put $R=R_\lambda$ and use the symplectic variables
\[
X=Rx,
\qquad
\Xi=R^{-1}\xi.
\]
After removing the scalar order $R^{2m_2}$ from $p_\lambda$ and the
corresponding scalar factors from $w_+$ and $w_-$, the rescaled localized
symbol has uniformly bounded finite derivatives.  Moreover, on its support
$|\Xi|\asymp R$, whereas every nonzero derivative of either rescaled
Sobolev weight contributes at least one factor $R^{-1}$.  Derivatives of
$\psi(X/R)$ and $\psi'(X/R)$ also contribute factors $R^{-1}$.

Apply the oscillatory integral formula for the Weyl product successively
to the five factors.  Repeated use of
\[
(1-\Delta_U)e^{-2i\sigma(U,V)}
=(1+4|V|^2)e^{-2i\sigma(U,V)}
\]
and its analogue with $U$ and $V$ interchanged makes the oscillatory
integrals absolutely convergent after finitely many integrations by parts.
If the resulting point $(X,\Xi)$ is separated from the rescaled support of
$p_\lambda$, the same integrations by parts yield an arbitrary negative
power of that separation.  Differentiating the product only increases the
finite number of derivatives required from the factors.  Returning to
$(x,\xi)$ gives \eqref{eq:mesoscopic-conjugated-symbol}.  In particular,
only finitely many seminorms of $a$ and finitely many derivatives of
$T^{\pm1}$ on $K_1$ enter for each prescribed $A,B,L$.

For the kernel representation, set
\[
m:=\frac{x+y}{2},
\qquad
\theta_\lambda(m):=S_m\zeta_\lambda.
\]
In the Weyl kernel formula make the exact change of variables
\[
\xi=\theta_\lambda(m)+R_\lambda\eta.
\]
This gives \eqref{eq:mesoscopic-kernel-form}, where
\[
A_\lambda(m,u)
:=
(2\pi)^{-n}
\int_{\mathbb{R}^n}
 e^{iu\cdot\eta}
 b_\lambda
 \bigl(m,\theta_\lambda(m)+R_\lambda\eta\bigr)
\,\mathrm{d}\eta.
\]
Since $|\zeta_\lambda|\asymp R_\lambda^2$ and the derivatives of
$S_m$ are bounded on $K_2$,
\[
|\partial_m^\gamma\theta_\lambda(m)|
\leq
C_\gamma R_\lambda^2.
\]
Combining this with
\eqref{eq:mesoscopic-conjugated-symbol}, each $m$-derivative of the
centered symbol costs at most one factor $R_\lambda$.  Integrating by parts
in $\eta$ yields \eqref{eq:mesoscopic-amplitude-decay}.  The spatial tails
outside $K_2$ follow from the phase-space decay already obtained in
\eqref{eq:mesoscopic-conjugated-symbol}.
\end{proof}

\begin{proposition}[Uniform rescaling of the Sobolev-conjugated block]
\label{prop:weyl-local-uniform-semiclass}
Let $\widetilde T_{j,k}$ be as in
Definition~\ref{def:sobolev-conjugated-blocks}.  Put $R:=R_k$ and define
the unitary dilation
\[
(U_Ru)(X):=R^{-n/2}u(X/R).
\]
If $b_{j,k}$ is the exact Weyl symbol from
Lemma~\ref{lem:mesoscopic-kernel-localization}, then
\[
U_R\widetilde T_{j,k}U_R^{-1}
=
\operatorname{Op}^w(B_{j,k}),
\qquad
B_{j,k}(X,\Xi):=b_{j,k}(X/R,R\Xi).
\]
For every finite derivative order required by Calder\'on--Vaillancourt,
there are finite integers $A,B$ such that
\[
\sup_{j,k}\sup_{X,\Xi}
|\partial_X^\alpha\partial_\Xi^\beta B_{j,k}(X,\Xi)|
\leq
C_{\alpha,\beta}\mathcal S_{A,B}(a;K_1).
\]
Consequently,
\begin{equation}\label{eq:semi-main-updated}
\|\widetilde T_{j,k}\|_{L^2\to L^2}
\leq
C\mathcal S_{A,B}(a;K_1)
\end{equation}
uniformly in $j,k$.  Equivalently,
\[
\|T_{j,k}\|_{H^s\to H^{s-m_2}}
\leq
C\mathcal S_{A,B}(a;K_1).
\]
If the target Sobolev order is not compensated by $m_2$, then
\begin{equation}\label{eq:semi-same-order-band}
\|T_{j,k}\|_{H^s\to H^s}
\leq
C R_k^{2m_2}\mathcal S_{A,B}(a;K_1)
=
C2^{km_2}\mathcal S_{A,B}(a;K_1).
\end{equation}
\end{proposition}

\begin{proof}
The linear map
\[
(x,\xi)=(X/R,R\Xi)
\]
is symplectic, so Weyl quantization is exactly covariant under $U_R$.
Taking $L=0$ in
\eqref{eq:mesoscopic-conjugated-symbol} gives
\[
|\partial_x^\alpha\partial_\xi^\beta b_{j,k}(x,\xi)|
\leq
C_{\alpha,\beta}
R^{|\alpha|-|\beta|}
\mathcal S_{A,B}(a;K_1)
\]
at every prescribed finite derivative level.  Therefore
\[
|\partial_X^\alpha\partial_\Xi^\beta B_{j,k}(X,\Xi)|
=
R^{-|\alpha|+|\beta|}
|\partial_x^\alpha\partial_\xi^\beta b_{j,k}(X/R,R\Xi)|
\leq
C_{\alpha,\beta}\mathcal S_{A,B}(a;K_1).
\]
Calder\'on--Vaillancourt gives \eqref{eq:semi-main-updated}.  The
same-order estimate follows from
$\langle\xi\rangle\asymp R_k^2$ on the block.
\end{proof}

\begin{proposition}[Uniform local Weyl estimate on a mesoscopic block]
\label{prop:weyl-local-uniform-latex:updated}
Let $a\in S^{m_1,m_2}$, let $K\subset\mathbb{R}^n$ be compact, and let
$\chi,\chi'\in C_c^\infty(K)$.  Fix $s\in\mathbb{R}$.  Then there exist
integers $A=A(n,s,m_2)$ and $B=B(n,s,m_2)$ and a compact
$K_1\supset K$ such that
\begin{equation}\label{eq:weyl-local-upd}
\|\chi\operatorname{Op}^w(a\Lambda_{j,k})\chi'\|_{H^s\to H^{s-m_2}}
\leq
C\mathcal S_{A,B}(a;K_1)
\end{equation}
for all $k\geq k_\ast$, uniformly in $j,k$.  The fixed low-frequency
block and the finitely many blocks with $0\leq k<k_\ast$ satisfy the same
mapping estimate after enlarging the constant.
\end{proposition}

\begin{proof}
Apply Proposition~\ref{prop:weyl-local-uniform-semiclass} with
$\psi=\chi$ and $\psi'=\chi'$.  This gives the high-frequency estimate.
The low-frequency and transition blocks are treated by
Lemma~\ref{lem:local-sobolev-cv-reduction}.
\end{proof}

\begin{remark}[No positive band loss at the natural Sobolev order]
The estimate \eqref{eq:weyl-local-upd} is uniform in the block index.  In
particular, the factor $2^{kN_\xi}$ occurring for the former fixed-radius
partition is absent.  The cancellation is scale-theoretic: the block
varies at the phase-space scale
\[
R_k^2|\mathrm{d}x|^2+R_k^{-2}|\mathrm{d}\xi|^2,
\]
whereas the Sobolev weights vary in frequency at the larger dyadic scale
$R_k^2$.
\end{remark}

\begin{corollary}[Order-zero case]
\label{cor:weyl-m2zero:conservative}
If $m_2=0$, then
\[
\|\chi\operatorname{Op}^w(a\Lambda_{j,k})\chi'\|_{H^s\to H^s}
\leq
C\mathcal S_{A,B}(a;K_1)
\]
uniformly in $j$ and $k\geq0$.
\end{corollary}

\begin{proof}
This is Proposition~\ref{prop:weyl-local-uniform-latex:updated} with
$m_2=0$.
\end{proof}

\subsection{Global weighted Sobolev bounds}
\label{subsec:global-weighted-sobolev}

\begin{definition}[Weighted Sobolev spaces]
\label{def:weighted-sobolev-spaces}
For $s,\sigma\in\mathbb{R}$ define
\[
H^{s,\sigma}(\mathbb{R}^n)
:=
\left\{
u\in\mathcal S'(\mathbb{R}^n):
\langle x\rangle^\sigma\langle D\rangle^s u\in L^2(\mathbb{R}^n)
\right\},
\]
with norm
\[
\|u\|_{H^{s,\sigma}}
:=
\|\langle x\rangle^\sigma\langle D\rangle^s u\|_{L^2}.
\]
\end{definition}

\begin{theorem}[Global weighted Sobolev estimate]
\label{thm:global-weighted-sobolev}
Let $a\in S^{m_1,m_2}$ and let $s,\sigma\in\mathbb{R}$.  Then
$\operatorname{Op}^w(a)$ extends uniquely to a bounded operator
\[
\operatorname{Op}^w(a):
H^{s,\sigma}(\mathbb{R}^n)
\longrightarrow
H^{s-m_2,\sigma-m_1}(\mathbb{R}^n).
\]
More precisely,
\begin{equation}\label{eq:global-weighted-sobolev}
\|\operatorname{Op}^w(a)u\|_{H^{s-m_2,\sigma-m_1}}
\leq
C_{n,s,\sigma,m_1,m_2}
\mathcal S_{D_{\operatorname{wSob}},D_{\operatorname{wSob}}}(a)
\|u\|_{H^{s,\sigma}},
\end{equation}
where the explicit sufficient derivative budget
$D_{\operatorname{wSob}}$ is defined in
\eqref{eq:weighted-sobolev-budget}.
\end{theorem}

\begin{proof}
For $u\in\mathcal S(\mathbb{R}^n)$ set
\[
v
:=
\langle x\rangle^\sigma\langle D\rangle^s u.
\]
Then
\[
u
=
\langle D\rangle^{-s}\langle x\rangle^{-\sigma}v.
\]
Hence
\[
\|\operatorname{Op}^w(a)u\|_{H^{s-m_2,\sigma-m_1}}
=
\|\mathcal A_{s,\sigma}v\|_{L^2},
\]
where
\[
\mathcal A_{s,\sigma}
:=
\langle x\rangle^{\sigma-m_1}
\langle D\rangle^{s-m_2}
\operatorname{Op}^w(a)
\langle D\rangle^{-s}
\langle x\rangle^{-\sigma}.
\]
Since multiplication by a function of $x$ and Fourier multiplication by a
function of $\xi$ are Weyl operators with the corresponding symbols, the
exact Weyl symbol of $\mathcal A_{s,\sigma}$ is
\[
c_{s,\sigma}
=
\langle x\rangle^{\sigma-m_1}
\#
\langle\xi\rangle^{s-m_2}
\# a\#
\langle\xi\rangle^{-s}
\#
\langle x\rangle^{-\sigma}.
\]
The spatial orders cancel,
\[
(\sigma-m_1)+m_1-\sigma=0,
\]
and the frequency orders cancel,
\[
(s-m_2)+m_2-s=0.
\]
Proposition~\ref{prop:weighted-sobolev-derivative-budget} therefore gives
$c_{s,\sigma}\in S^{0,0}$ and controls all derivatives required by
Theorem~\ref{thm:CV-semiclass-uniform} by
\[
C
\mathcal S_{D_{\operatorname{wSob}},D_{\operatorname{wSob}}}(a).
\]
Applying Theorem~\ref{thm:CV-semiclass-uniform} with $h=1$ yields
\[
\|\mathcal A_{s,\sigma}\|_{L^2\to L^2}
\leq
C
\mathcal S_{D_{\operatorname{wSob}},D_{\operatorname{wSob}}}(a).
\]
This proves \eqref{eq:global-weighted-sobolev} for Schwartz functions, and
the extension follows by density.
\end{proof}

\begin{remark}[Global operator bound versus blockwise mesoscopic bounds]
\label{rem:weighted-global-vs-blocks}
Theorem~\ref{thm:global-weighted-sobolev} concerns the complete symbol
$a\in S^{m_1,m_2}$ and does not require spatial compactness or any
additional hypothesis on $G_x$.  It does not imply a uniform natural-order
estimate
\[
\operatorname{Op}^w(a\Lambda_{j,k}):
H^{s,\sigma}
\longrightarrow
H^{s-m_2,\sigma-m_1}
\]
for the individual mesoscopic blocks under the present assumptions.  Indeed,
Proposition~\ref{prop:Lambda-deriv-bounds:full} gives globally
\[
|\partial_x^\alpha\partial_\xi^\beta\Lambda_{j,k}(x,\xi)|
\leq
C_{\alpha,\beta}
\langle x\rangle^{|\alpha|}
R_k^{|\alpha|-|\beta|},
\]
so the blockwise global $x$-derivatives need not preserve the spatial order
$m_1$.  The natural-order mesoscopic estimates used for almost orthogonality
therefore remain spatially localized, whereas the weighted Sobolev estimate
above is a global property of the full product-type symbol class
$S^{m_1,m_2}$.
\end{remark}

\subsection{Composition and recombination}
\label{subsec:moyal-truncation-dyadic}

\begin{lemma}[Off-diagonal decay and Cotlar--Stein summability]
\label{lem:mesoscopic-offdiagonal-cotlar}
Under the hypotheses of
Lemma~\ref{lem:mesoscopic-kernel-localization}, for every
$M\in\mathbb{N}_0$ there is a constant $C_M$, depending on finitely many
seminorms of $a$ and finitely many derivatives of $T^{\pm1}$ on a fixed
compact set, such that
\begin{equation}\label{eq:mesoscopic-offdiag-adjoint}
\|\widetilde T_\lambda^*\widetilde T_\mu\|_{L^2\to L^2}
\leq
C_M
\left(
1+\operatorname{dist}_{\mathrm{mes}}(\lambda,\mu)
\right)^{-M}
\end{equation}
and
\begin{equation}\label{eq:mesoscopic-offdiag-dual}
\|\widetilde T_\lambda\widetilde T_\mu^*\|_{L^2\to L^2}
\leq
C_M
\left(
1+\operatorname{dist}_{\mathrm{mes}}(\lambda,\mu)
\right)^{-M}.
\end{equation}
For the quantitative Cotlar--Stein choice
\[
M=M_{\operatorname{CS}}:=2n+2,
\]
all constants in the preceding estimates are controlled by
\[
\mathcal S_{D_{\operatorname{CS}},D_{\operatorname{CS}}}(a;K_1)
\]
and by derivatives of $T$ and $T^{-1}$ through order
$D_{\operatorname{CS}}$, where
\[
D_{\operatorname{CS}}
=
2D_{\operatorname{ker}}
+
8Q_{\operatorname{osc}}
+
M_{\operatorname{CS}}
+
2
\]
is the explicit sufficient derivative budget defined in
\eqref{eq:CS-derivative-budget}.  See
Proposition~\ref{prop:mesoscopic-cotlar-derivative-budget}.
Moreover, for every $P>n$,
\begin{equation}\label{eq:mesoscopic-discrete-sum}
\sup_\lambda
\sum_\mu
\left(
1+\operatorname{dist}_{\mathrm{mes}}(\lambda,\mu)
\right)^{-P}
<\infty,
\end{equation}
and the analogous estimate with $\lambda$ and $\mu$ interchanged holds.
Consequently, every integer $M>2n$ yields the two Cotlar--Stein sums
\begin{equation}\label{eq:cotlar-hyp-1}
\sup_\lambda
\sum_\mu
\|\widetilde T_\lambda^*\widetilde T_\mu\|_{L^2\to L^2}^{1/2}
<\infty
\end{equation}
and
\begin{equation}\label{eq:cotlar-hyp-2}
\sup_\mu
\sum_\lambda
\|\widetilde T_\lambda\widetilde T_\mu^*\|_{L^2\to L^2}^{1/2}
<\infty.
\end{equation}
\end{lemma}

\begin{proof}
We prove \eqref{eq:mesoscopic-offdiag-adjoint}; the second estimate is
symmetric.  The kernel of $\widetilde T_\lambda^*\widetilde T_\mu$ is an
integral in $z$ with phase
\[
\Phi(z;x,y)
:=
(z-y)\cdot S_{(z+y)/2}\zeta_\mu
-
(z-x)\cdot S_{(z+x)/2}\zeta_\lambda.
\]
Set
\[
\Delta\zeta:=\zeta_\mu-\zeta_\lambda,
\qquad
\Delta:=|\Delta\zeta|,
\qquad
v:=\frac{\Delta\zeta}{\Delta}
\]
when $\Delta\neq0$.  If
$\operatorname{dist}_{\mathrm{mes}}(\lambda,\mu)$ is bounded by a fixed
structural constant, the estimate follows from the uniform block bound
\eqref{eq:semi-main-updated}.  We therefore assume that this distance is
larger than that constant.

For
\[
H(z,x):=S_{(z+x)/2}(z-x)
\]
one has, on a fixed compact neighborhood of the spatial supports,
\[
D_zH(z,x)=S_z+E(z,x),
\qquad
\|E(z,x)\|\leq C|z-x|.
\]
Hence
\[
\nabla_z\Phi
=
S_z\Delta\zeta
+
E(z,y)^\top\zeta_\mu
-
E(z,x)^\top\zeta_\lambda.
\]
Let
\[
d:=\operatorname{dist}_{\mathrm{mes}}(\lambda,\mu)
=
\frac{\Delta}{R_\lambda+R_\mu}.
\]
In the effective interaction region
\[
R_\lambda|z-x|\leq\varepsilon d,
\qquad
R_\mu|z-y|\leq\varepsilon d,
\]
we use $|\zeta_\lambda|\asymp R_\lambda^2$ and
$|\zeta_\mu|\asymp R_\mu^2$ to obtain
\[
|E(z,y)^\top\zeta_\mu-E(z,x)^\top\zeta_\lambda|
\leq
C\varepsilon d(R_\lambda+R_\mu)
=
C\varepsilon\Delta.
\]
Uniform ellipticity of $S_z$ gives a constant $c_S>0$ such that
\[
v^\top S_zv\geq c_S.
\]
Choosing $\varepsilon>0$ sufficiently small therefore gives
\begin{equation}\label{eq:phase-gradient-lower}
|v\cdot\nabla_z\Phi|
\geq
\frac{c_S}{2}\Delta.
\end{equation}
For every prescribed finite derivative order, Taylor expansion of
$D_zH(z,x)$ at $x=z$ shows that the corresponding derivatives of
$v\cdot\nabla_z\Phi$ are bounded by $C\Delta$ in the same region: the
terms independent of $z-x$ and $z-y$ multiply $\Delta\zeta$, whereas all
remaining terms retain one factor $|z-x|$ or $|z-y|$ and are controlled
by the preceding estimate.  Consequently the derivatives of
$(v\cdot\nabla_z\Phi)^{-1}$ required below are bounded by $C\Delta^{-1}$.

Introduce
\[
\mathcal L
:=
\frac{1}{i\,v\cdot\nabla_z\Phi}
\,v\cdot\nabla_z,
\qquad
\mathcal L e^{i\Phi}=e^{i\Phi}.
\]
By Lemma~\ref{lem:mesoscopic-kernel-localization}, every $z$-derivative of
the product amplitude costs at most $C(R_\lambda+R_\mu)$.  Thus each
integration by parts with $\mathcal L$ contributes a factor
\[
\frac{R_\lambda+R_\mu}{\Delta}=d^{-1}.
\]
After $M$ integrations by parts, the part of the kernel in the effective
interaction region is bounded by $C_Md^{-M}$ times an integrable
convolution of the two kernel envelopes.

Outside that region, at least one of
$R_\lambda|z-x|$ and $R_\mu|z-y|$ is larger than $\varepsilon d$.
Choosing the exponent $N$ in
\eqref{eq:mesoscopic-amplitude-decay} larger than $M+n+1$ extracts the
same factor $(1+d)^{-M}$ while leaving an integrable kernel envelope.
The spatial tails outside the fixed compact neighborhood are treated by
the arbitrary phase-space decay in
\eqref{eq:mesoscopic-conjugated-symbol}.  Therefore, for some $N_0>n$,
\[
|K_{\lambda,\mu}(x,y)|
\leq
C_M(1+d)^{-M}
R_\lambda^nR_\mu^n
\int
\langle R_\lambda(z-x)\rangle^{-N_0}
\langle R_\mu(z-y)\rangle^{-N_0}
\,\mathrm{d}z.
\]
Integrating first in $y$ and then in $z$, and symmetrically first in $x$,
gives
\[
\sup_x\int|K_{\lambda,\mu}(x,y)|\,\mathrm{d}y
+
\sup_y\int|K_{\lambda,\mu}(x,y)|\,\mathrm{d}x
\leq
C_M(1+d)^{-M}.
\]
Schur's test proves \eqref{eq:mesoscopic-offdiag-adjoint}.

It remains to verify the discrete sum.  If $R_\mu\asymp R_\lambda$, only
finitely many dyadic values of $k'$ occur.  The centers at each such scale
are uniformly separated at distance comparable to $R_\lambda$, so a
packing argument gives
\[
\#\{\mu:q\leq d(\lambda,\mu)<q+1,\ R_\mu\asymp R_\lambda\}
\leq
C(1+q)^{n-1}.
\]
Hence the comparable-scale contribution to
\eqref{eq:mesoscopic-discrete-sum} converges for $P>n$.

If $R_\mu\geq C_0R_\lambda$ for a sufficiently large structural
constant, radial dyadic separation gives
\[
|\zeta_\mu-\zeta_\lambda|\asymp R_\mu^2,
\qquad
d(\lambda,\mu)\asymp R_\mu.
\]
Since Section~\ref{subsec:dyadic-decomposition} gives
$\#J_{k'}\leq CR_\mu^n$, the contribution of the $k'$-th scale is bounded
by $CR_\mu^{n-P}$ and is geometrically summable for $P>n$.  The case
$R_\mu\leq C_0^{-1}R_\lambda$ is symmetric.  This proves
\eqref{eq:mesoscopic-discrete-sum}.  Taking $P=M/2$ gives
\eqref{eq:cotlar-hyp-1}--\eqref{eq:cotlar-hyp-2} whenever $M>2n$.
\end{proof}

\begin{proposition}[Quantitative Cotlar--Stein recombination]
\label{thm:global-sum-semiclassical}
Let $a\in S^{m_1,m_2}$, fix $s\in\mathbb{R}$, and let
$\psi,\psi'\in C_c^\infty(\mathbb{R}^n)$.  For the mesoscopic partition
constructed in Section~\ref{subsec:dyadic-decomposition}, the
high-frequency blocks
\[
T_{j,k}
=
\psi\operatorname{Op}^w(a\Lambda_{j,k})\psi'
\]
satisfy the Cotlar--Stein conditions
\eqref{eq:cotlar-hyp-1}--\eqref{eq:cotlar-hyp-2}.  In particular, one may
fix
\[
M_{\operatorname{CS}}:=2n+2.
\]
Then
\[
\sum_{k\geq0}\sum_{j\in J_k}T_{j,k}
\]
converges strongly from $H^s$ to $H^{s-m_2}$ and defines a bounded
operator
\[
S:H^s\to H^{s-m_2}.
\]
Its norm is controlled by finitely many local seminorms of $a$, the
ellipticity constants, finitely many derivatives of $T^{\pm1}$ on a fixed
compact neighborhood of the spatial cutoffs, and finitely many seminorms
of $\psi,\psi'$.

Moreover,
\[
S+
\psi\operatorname{Op}^w(a\Lambda_0)\psi'
=
\psi\operatorname{Op}^w(a)\psi'
\]
in $\mathcal L(H^s,H^{s-m_2})$.
\end{proposition}

\begin{proof}
Apply Lemma~\ref{lem:mesoscopic-offdiagonal-cotlar} with
$M=M_{\operatorname{CS}}=2n+2$.  Then
$M/2=n+1>n$, so
\eqref{eq:cotlar-hyp-1}--\eqref{eq:cotlar-hyp-2} hold.  The
Cotlar--Stein lemma, Lemma~\ref{lem:cotlar-stein}, gives boundedness and
strong convergence of
\[
\sum_{k\geq0}\sum_{j\in J_k}\widetilde T_{j,k}
\]
on $L^2$.  Conjugating back gives boundedness and strong convergence of
the corresponding series from $H^s$ to $H^{s-m_2}$.

Finally,
\[
\Lambda_0+
\sum_{k\geq0}\sum_{j\in J_k}\Lambda_{j,k}=1
\]
by \eqref{eq:dyadic-partition-unity}.  On Schwartz functions, local
finiteness of the partition in phase space and the standard oscillatory
interpretation of Weyl quantization identify the strong limit with
$\psi\operatorname{Op}^w(a(1-\Lambda_0))\psi'$.  Both sides extend
boundedly from $H^s$ to $H^{s-m_2}$, which proves the stated identity.
\end{proof}

\begin{remark}[Nature of the off-diagonal attenuation]
The interaction is not exactly orthogonal for separated blocks.  The
quantitative conclusion is rapid decay in the normalized mesoscopic
distance \eqref{eq:mesoscopic-block-distance}.  The scale
$R_k=2^{k/2}$ is essential: the kernel of a block is localized at spatial
scale $R_k^{-1}$, and the variation of the physical center
$T_x^{-1}\zeta_{j,k}$ across that scale is of order $R_k$, exactly the
frequency width of the block.  This balance is what makes the lower bound
\eqref{eq:phase-gradient-lower} uniform in $j$ and $k$.
\end{remark}

%% file: 50_app1_dyadic.tex
\section{Model uses of the localized calculus}
\label{sec:applications}

The purpose of this section is not to claim new global anisotropic parametrices beyond the finite-seminorm calculus proved above.  The results below record two model uses of the mesoscopic localization: first, a microlocal parametrix scheme for elliptic Weyl symbols; second, a compatibility statement for the Euclidean Radon transform as a model Fourier integral operator.  The global recombination of the pseudodifferential blocks is no longer imposed as an additional hypothesis: it follows from the quantitative off-diagonal estimates and Cotlar--Stein summability proved in Section~\ref{sec:quantization-local-global}.  The Radon recombination is then obtained by composition with a fixed localized Fourier integral operator.

Throughout this section $T_x=G_x^{1/2}$ and
\[
\|\xi\|_{g_x}=|T_x\xi|.
\]
By Proposition~\ref{prop:gx-basic}, the uniform ellipticity of $G_x$ gives
\[
1+\|\xi\|_{g_x}\asymp \langle\xi\rangle,
\]
with constants depending only on the ellipticity constants of $G_x$.  Thus the anisotropy used here is geometric and patchwise: it enters through the moving normalized variable $T_x\xi$ and through the mesoscopic derivative balance of the cutoffs $\Lambda_{j,k}$.

The dyadic partition is the one constructed in Section~\ref{subsec:dyadic-decomposition}.  It has the form
\[
\Lambda_0+
\sum_{k\geq0}\sum_{j\in J_k}\Lambda_{j,k}=1,
\]
where $\Lambda_0$ is a fixed low-frequency block and the indices $k\geq0$ describe the high-frequency dyadic part.  The estimates below are stated for the high-frequency blocks.  The low-frequency block is always treated separately by the standard compact-frequency pseudodifferential calculus.

\subsection{Microlocal parametrix for symbols with \texorpdfstring{$x$}{x}-dependent fiber scale}
\label{sec:parametrix}

Fix a compact set $K\subset\mathbb{R}^n_x$.  Choose
\[
\psi,\psi',\chi\in C_c^\infty(K)
\]
so that $\chi\equiv1$ on a neighborhood of
$\operatorname{supp}\psi\cup\operatorname{supp}\psi'$.
Let $p\in S^{m_1,m_2}$, with $m_2>0$, and assume that $p$ is microlocally elliptic on a neighborhood of these supports: there exist constants $c_e>0$ and $R_e\geq1$ such that
\[
|p(x,\xi)|\geq c_e\langle\xi\rangle^{m_2},
\qquad
x\in K,
\quad
\langle\xi\rangle\geq R_e.
\]
Equivalently, after changing constants, one may write the lower bound with $1+\|\xi\|_{g_x}$ in place of $\langle\xi\rangle$.

By the standard elliptic symbolic construction in the Weyl calculus, choose a symbol
\[
q\in S^{-m_1,-m_2}
\]
which is a microlocal inverse for $p$ on a neighborhood of
$\operatorname{supp}\psi\cup\operatorname{supp}\psi'$, so that
\[
q\# p=1+r_\infty
\]
there, with $r_\infty$ smoothing; see, for example, \cite{HormanderIII}.  For $k\geq0$ define
\[
q_{j,k}:=q\Lambda_{j,k},
\qquad
Q_{j,k}:=
\psi\operatorname{Op}^w(q_{j,k})\chi,
\]
and set
\[
Q_0:=\psi\operatorname{Op}^w(q\Lambda_0)\chi.
\]

\begin{theorem}[Localized mesoscopic parametrix]
\label{thm:parametrix-micro}
Let $s\in\mathbb{R}$ and assume the hypotheses above.  Then the high-frequency series
\[
\sum_{k\geq0}\sum_{j\in J_k}Q_{j,k}
\]
converges strongly from $H^{s-m_2}$ to $H^s$.  Moreover,
\[
Q:=Q_0+
\sum_{k\geq0}\sum_{j\in J_k}Q_{j,k}
=
\psi\operatorname{Op}^w(q)\chi
\]
in $\mathcal L(H^{s-m_2},H^s)$.

If
\[
P_{\chi,\psi'}:=
\chi\operatorname{Op}^w(p)\psi',
\]
then there exists a smoothing operator $R$ such that
\[
Q P_{\chi,\psi'}
=
\psi\psi'+R
\]
as operators $H^s\to H^s$.  In particular,
\[
Q:H^{s-m_2}\longrightarrow H^s
\]
is a localized parametrix of order $-m_2$ modulo smoothing terms.

For suitable finite integers $A,B$ and a fixed compact neighborhood $K_1\supset K$,
\[
\|Q\|_{H^{s-m_2}\to H^s}
\leq
C\,\mathcal S_{A,B}(q;K_1),
\]
where $C$ depends only on the structural parameters, the spatial cutoffs, the ellipticity constants of $G_x$, and finitely many derivatives of $T_x^{\pm1}$ on $K_1$.  At every fixed finite derivative level, the standard elliptic construction controls the required seminorms of $q$ by finitely many local seminorms of $p$ and the ellipticity lower bound.
\end{theorem}

\begin{proof}
The symbol $q$ has fiber order $-m_2$.  Apply Theorem~\ref{thm:global-sum-semiclassical} to $a=q$, with source Sobolev index $s-m_2$ and spatial cutoffs $\psi,\chi$.  Since
\[
(s-m_2)-(-m_2)=s,
\]
the high-frequency series converges strongly from $H^{s-m_2}$ to $H^s$, and after adding the low-frequency block one obtains
\[
Q_0+
\sum_{k\geq0}\sum_{j\in J_k}Q_{j,k}
=
\psi\operatorname{Op}^w(q)\chi.
\]
No additional Cotlar--Stein hypothesis is required: the off-diagonal attenuation and the discrete summability are already part of Theorem~\ref{thm:global-sum-semiclassical}.

For comparison with the finite-seminorm formulation, Proposition~\ref{prop:mult-by-Lambda-preserves-S-anisotropic:full} gives, at derivative level $(\imath,\ell)$,
\[
N_x=2\imath,
\qquad
N_\xi=\left\lceil\frac{\imath+\ell}{2}\right\rceil.
\]
These are conservative global finite-level bookkeeping losses.  The parametrix mapping itself is governed by the sharper mesoscopic block estimate of Proposition~\ref{prop:weyl-local-uniform-latex:updated}, so no positive dyadic band loss occurs at the natural Sobolev order.

It remains to identify the localized composition.  By definition,
\[
Q P_{\chi,\psi'}
=
\psi\operatorname{Op}^w(q)\chi^2
\operatorname{Op}^w(p)\psi'.
\]
Because $\chi\equiv1$ on a neighborhood of $\operatorname{supp}\psi'$, the operator
\[
(1-\chi^2)\operatorname{Op}^w(p)\psi'
\]
is smoothing: its left and right spatial supports are separated, so the pseudodifferential kernel is evaluated away from the diagonal.  Hence
\[
Q P_{\chi,\psi'}
=
\psi\operatorname{Op}^w(q)\operatorname{Op}^w(p)\psi'
+R_1
\]
with $R_1$ smoothing.  The Weyl composition formula and the microlocal identity $q\# p=1+r_\infty$ then give
\[
\psi\operatorname{Op}^w(q)\operatorname{Op}^w(p)\psi'
=
\psi\operatorname{Op}^w(q\# p)\psi'
=
\psi\psi'+R_2,
\]
where $R_2$ is smoothing after the fixed spatial cutoffs are inserted.  Taking $R=R_1+R_2$ proves the parametrix identity.
\end{proof}

\begin{remark}[Direction of the Sobolev mapping]
The direction $Q:H^{s-m_2}\to H^s$ is essential.  The operator $\operatorname{Op}^w(p)$ has order $m_2$ in the fiber variable and therefore maps $H^s$ to $H^{s-m_2}$ locally.  A parametrix has the opposite order and acts in the reverse Sobolev direction.  In the recombination theorem this is exactly the specialization of the general mapping law to a symbol of fiber order $-m_2$.
\end{remark}

\begin{remark}[Dependence of constants]
The constants above are finite-level constants.  They may depend on $K$, $s$, $n$, $m_1$, $m_2$, the cutoffs $\psi,\psi',\chi$, the dyadic cutoff profile, the overlap and ellipticity constants, finitely many derivatives of $T_x^{\pm1}$ on a compact neighborhood, and finitely many local seminorms of $q$.  Through the fixed finite-order elliptic construction, the latter depend on finitely many local seminorms of $p$ and on the ellipticity lower bound.  No independent Cotlar--Stein constants enter the statement.
\end{remark}

\begin{example}[Local Schr\"odinger model]
\label{ex:schrodinger-local-new}
Let $K\subset\mathbb{R}^n$ be compact and let $\lambda\in C^\infty(\mathbb{R}^n)$ have bounded derivatives on $K$.  For
\[
p(x,\xi)=|\xi|^2+\lambda(x)
\]
with $G_x=\operatorname{Id}$, the preceding theorem applies on any microlocally elliptic region.  It gives a localized parametrix
\[
Q:H^{s-2}\to H^s
\]
modulo smoothing remainders.  The block recombination is supplied by Theorem~\ref{thm:global-sum-semiclassical}; no separate almost-orthogonality hypothesis is imposed in this example.
\end{example}

\begin{example}[Quadratic symbol with variable metric]
\label{ex:anisotropic-quadratic-new}
Let $K\subset\mathbb{R}^n$ be compact.  Assume $B\in C^\infty(\mathbb{R}^n;\operatorname{Sym}^+(n))$ satisfies
\[
c|\eta|^2\leq \eta^\top B(x)\eta\leq C|\eta|^2,
\qquad x\in K,
\]
and has bounded derivatives on $K$.  Let $\lambda\in C^\infty(\mathbb{R}^n)$ have bounded derivatives on $K$, and set
\[
p(x,\xi)=\xi^\top B(x)\xi+\lambda(x).
\]
With $G_x=B(x)$, the fiber norm is $\|\xi\|_{g_x}=(\xi^\top B(x)\xi)^{1/2}$ and remains uniformly comparable with $|\xi|$ on $K$.  On a microlocally elliptic region, the theorem gives a localized parametrix $Q:H^{s-2}\to H^s$ modulo smoothing terms.  This example illustrates the patchwise moving geometry of the construction, not a new global order scale distinct from the Euclidean one.
\end{example}

\subsection{Radon transform and dyadic microlocal decomposition}
\label{sec:radon-dyadic}

Let $R$ denote the Euclidean Radon transform
\[
(Rf)(s,\omega)
=
\int_{\mathbb{R}^n}f(x)\delta(s-x\cdot\omega)\,\mathrm{d}x,
\qquad
(s,\omega)\in\mathbb{R}\times\mathbb{S}^{n-1}.
\]
As a Fourier integral operator, $R$ has order $-(n-1)/2$ and is associated with the standard Radon canonical relation; see \cite[Ch.~XI, Sections~2--3]{HormanderIV}.  The normal operator $R^*R$ is an elliptic pseudodifferential operator of order $-(n-1)$, with principal symbol proportional to $|\xi|^{-(n-1)}$ away from the zero section; see also \cite[Ch.~I--II]{HelgasonRadon}.

Set
\[
m_R:=\frac{n-1}{2}.
\]
Fix compact sets $K_x\subset\mathbb{R}^n$ and $K_{s,\omega}\subset\mathbb{R}\times\mathbb{S}^{n-1}$, and choose cutoffs
\[
\psi_x,\psi_x'\in C_c^\infty(K_x),
\qquad
\psi_{s,\omega}\in C_c^\infty(K_{s,\omega}).
\]
For $a\in S^{m_1,m_2}$ and $k\geq0$, define the localized Radon blocks
\[
T_{j,k}:=
\psi_{s,\omega}R\psi_x\operatorname{Op}^w(a\Lambda_{j,k})\psi_x'.
\]
The low-frequency block
\[
T_0:=
\psi_{s,\omega}R\psi_x\operatorname{Op}^w(a\Lambda_0)\psi_x'
\]
is treated separately by compact-frequency estimates.

\begin{proposition}[Uniform compensated local bound for Radon blocks]
\label{prop:radon-local-block}
For every $s\in\mathbb{R}$ there exist finite integers
$A=A(n,s,m_2)$ and $B=B(n,s,m_2)$, a compact neighborhood
$K_1\supset K_x$, and a constant
\[
C=C(K_x,K_{s,\omega},\psi_x,\psi_x',\psi_{s,\omega},n,s,m_1,m_2,\varphi,L)
\]
such that, for all $j$ and $k\geq0$,
\begin{equation}
\label{eq:radon-local-compensated}
\|T_{j,k}\|_{H^s_x\to H^{s+m_R-m_2}_{s,\omega}}
\leq
C\,\mathcal S_{A,B}(a;K_1).
\end{equation}
Thus the extra localization loss in the notation of the previous formulation is
\[
\boxed{N_{\operatorname{loc}}=0.}
\]
In particular,
\[
\begin{array}{c|c|c}
 n & m_R & \text{uniform compensated mapping} \\
\hline
 2 & \frac12 & H^s_x\to H^{s+\frac12-m_2}_{s,\omega} \\
 3 & 1 & H^s_x\to H^{s+1-m_2}_{s,\omega}
\end{array}
\]
and in both dimensions $N_{\operatorname{loc}}=0$.
\end{proposition}

\begin{proof}
Set
\[
A_{j,k}:=
\psi_x\operatorname{Op}^w(a\Lambda_{j,k})\psi_x'.
\]
Proposition~\ref{prop:weyl-local-uniform-latex:updated} gives
\[
\|A_{j,k}\|_{H^s_x\to H^{s-m_2}_x}
\leq
C\mathcal S_{A,B}(a;K_1)
\]
uniformly in $j$ and $k$.  Since $A_{j,k}u$ is supported in the fixed compact set determined by $\psi_x$, the compactly supported Fourier integral operator estimate for the Radon transform gives, with $t=s-m_2$,
\[
\|\psi_{s,\omega}Rv\|_{H^{t+m_R}_{s,\omega}}
\leq C_t\|v\|_{H^t_x}
\]
for $v$ supported in that compact set.  Therefore
\[
\|T_{j,k}u\|_{H^{s+m_R-m_2}_{s,\omega}}
\leq
C\mathcal S_{A,B}(a;K_1)\|u\|_{H^s_x},
\]
which proves \eqref{eq:radon-local-compensated}.  No positive power of $2^k$ is produced: the mesoscopic rescaling already absorbs the finite derivative costs of the block localization.  Hence $N_{\operatorname{loc}}=0$.
\end{proof}

\begin{proposition}[Non-compensated band estimate]
\label{prop:radon-local-band}
With the same hypotheses as Proposition~\ref{prop:radon-local-block}, one has
\begin{equation}
\label{eq:radon-local-band}
\|T_{j,k}\|_{H^s_x\to H^s_{s,\omega}}
\leq
C\,2^{k(m_2-m_R)}\mathcal S_{A,B}(a;K_1)
\end{equation}
for all $j$ and $k\geq0$.  Thus, explicitly,
\[
\begin{array}{c|c}
 n & \text{non-compensated band factor} \\
\hline
 2 & 2^{k(m_2-\frac12)} \\
 3 & 2^{k(m_2-1)}
\end{array}
\]
with no additional localization factor.
\end{proposition}

\begin{proof}
The Radon factor has order $-m_R$, whereas the pseudodifferential block has order $m_2$.  On the high-frequency block one has $\langle\xi\rangle\asymp2^k$.  The standard local FIO calculus therefore contributes the uncompensated factor $2^{k(m_2-m_R)}$.  Proposition~\ref{prop:radon-local-block} shows that there is no further positive factor coming from the mesoscopic localization.  The fixed low-frequency and transition blocks are absorbed into the constant.
\end{proof}

\begin{theorem}[Radon recombination]
\label{thm:radon-global}
Let $T_{j,k}$ be as above.  Then
\[
S_R:=
\sum_{k\geq0}\sum_{j\in J_k}T_{j,k}
\]
converges strongly from $H^s_x$ to
$H^{s+m_R-m_2}_{s,\omega}$ and defines a bounded operator satisfying
\begin{equation}
\label{eq:radon-global}
\|S_R\|_{H^s_x\to H^{s+m_R-m_2}_{s,\omega}}
\leq C\,\mathcal S_{D,D}(a;K_1)
\end{equation}
for some finite $D=D(n,s,m_2)$.
Moreover,
\[
S_R+T_0
=
\psi_{s,\omega}R\psi_x\operatorname{Op}^w(a)\psi_x'
\]
in
$\mathcal L(H^s_x,H^{s+m_R-m_2}_{s,\omega})$.
\end{theorem}

\begin{proof}
Apply Theorem~\ref{thm:global-sum-semiclassical} with
$\psi=\psi_x$ and $\psi'=\psi_x'$.  The series
\[
\sum_{k\geq0}\sum_{j\in J_k}
\psi_x\operatorname{Op}^w(a\Lambda_{j,k})\psi_x'
\]
converges strongly from $H^s_x$ to $H^{s-m_2}_x$ and, after adding the low-frequency block, equals
\[
\psi_x\operatorname{Op}^w(a)\psi_x'.
\]
The fixed localized Radon operator
\[
\psi_{s,\omega}R
\]
is bounded from $H^{s-m_2}_x$ to
$H^{s+m_R-m_2}_{s,\omega}$ on functions supported in the fixed compact set determined by $\psi_x$.  Composing this bounded operator with the strongly convergent series gives the asserted strong convergence and the identity after adding $T_0$.  No separate Cotlar--Stein hypothesis is required for the Radon family.
\end{proof}

\begin{remark}[Scope of the Radon model]
The Radon discussion is a compatibility model for the mesoscopic localization.  It uses the standard Euclidean FIO structure of $R$ and the uniformly comparable fiber scale $\|\xi\|_{g_x}\asymp |\xi|$.  The quantitative conclusion $N_{\operatorname{loc}}=0$ reflects the uniform mesoscopic block estimate of Section~\ref{sec:quantization-local-global}; it should not be read as a new anisotropic Radon calculus.  Any stronger statement would require additional hypotheses on the canonical relation or a genuinely non-uniform anisotropic geometry.
\end{remark}

%% file: 60_appendix_dyadic.tex
\appendix

\section{Appendix}\label{sec:appendix}

\subsection{Multivariable Fa\`a di Bruno lemma}

\begin{lemma}[Multivariable Fa\`a di Bruno: composition estimate]\label{lem:faa_bruno_estimate}
Let $m,n\in\mathbb{N}$. Let $\varphi\in C_c^\infty(\mathbb{R}^m)$ and let
$u\in C^\infty(\mathbb{R}^n_x\times\mathbb{R}^n_\xi;\mathbb{R}^m)$.
Fix a nonzero multi-index $\gamma=(\alpha,\beta)\in\mathbb{N}_0^{2n}\setminus\{0\}$, and set
$s:=|\gamma|=|\alpha|+|\beta|$. Assume that for every
$\gamma'=(\alpha',\beta')\in\mathbb{N}_0^{2n}$ there exists $C_{\gamma'}>0$ such that
\[
|\partial^{\gamma'}u(x,\xi)|
\leq C_{\gamma'}\,\langle x\rangle^{|\alpha'|}\,\langle\xi\rangle^{|\beta'|},
\qquad (x,\xi)\in\mathbb{R}^{2n}.
\]
Then
\[
\big|\partial^\gamma(\varphi\circ u)(x,\xi)\big|
\leq C_\gamma(\varphi)\,\langle x\rangle^{|\alpha|}\,\langle\xi\rangle^{|\beta|},
\qquad (x,\xi)\in\mathbb{R}^{2n},
\]
where
\[
C_\gamma(\varphi)
:=\sum_{r=1}^s\frac{p_r(\varphi)}{r!}
\sum_{(\gamma^1,\dots,\gamma^r)\in\mathcal{S}_r(\gamma)}
\prod_{\nu=1}^r C_{\gamma^\nu},
\]
\[
\mathcal{S}_r(\gamma)
:=\Big\{(\gamma^1,\dots,\gamma^r)\in(\mathbb{N}_0^{2n}\setminus\{0\})^r:
\sum_{\nu=1}^r\gamma^\nu=\gamma\Big\},
\]
and $p_r(\varphi):=\sup_{y\in\mathbb{R}^m}\|\nabla^r\varphi(y)\|_{\operatorname{op}}$.
Moreover, $\partial^\gamma(\varphi\circ u)\equiv0$ outside $u^{-1}(\operatorname{supp}\varphi)$.
The constant $C_\gamma(\varphi)$ depends only on $p_1(\varphi),\dots,p_s(\varphi)$ and on
$\{C_{\gamma'}:|\gamma'|\leq s\}$.
\end{lemma}

\begin{proof}
By the ordered multivariable Fa\`a di Bruno formula; see, for instance,
\cite{ConstantineSavitsFaaDiBruno},
\[
\partial^\gamma(\varphi\circ u)
=\sum_{r=1}^s\frac{1}{r!}
\sum_{\substack{\gamma^1+\cdots+\gamma^r=\gamma\\ \gamma^\nu\neq0}}
\nabla^r\varphi(u(x,\xi))
\big(\partial^{\gamma^1}u(x,\xi),\dots,\partial^{\gamma^r}u(x,\xi)\big).
\]
Taking absolute values, using the definition of $p_r(\varphi)$, and applying the hypotheses on the
partial derivatives of $u$ gives the asserted estimate. The support statement follows from
$\operatorname{supp}(\varphi\circ u)\subset u^{-1}(\operatorname{supp}\varphi)$.
\end{proof}

\subsection{Weyl composition in anisotropic classes and semiclassical version}
\label{sec:appendix-weyl-moyal}

For $m_1,m_2\in\mathbb{R}$ we denote by
$S^{m_1,m_2}=S^{m_1,m_2}(\mathbb{R}^n_x\times\mathbb{R}^n_\xi)$ the class of all
$a\in C^\infty(\mathbb{R}^{2n})$ such that, for every pair of multi-indices $\alpha,\beta$,
\begin{equation}\label{eq:anisotropic-class}
\big|\partial_x^\alpha\partial_\xi^\beta a(x,\xi)\big|
\leq C_{\alpha,\beta}\,\langle x\rangle^{m_1-|\alpha|}\,\langle\xi\rangle^{m_2-|\beta|},
\qquad (x,\xi)\in\mathbb{R}^{2n}.
\end{equation}
We use the seminorms
\[
\|a\|_{\alpha,\beta}^{(m_1,m_2)}
:=\sup_{(x,\xi)\in\mathbb{R}^{2n}}
\langle x\rangle^{-m_1+|\alpha|}\langle\xi\rangle^{-m_2+|\beta|}
\big|\partial_x^\alpha\partial_\xi^\beta a(x,\xi)\big|.
\]
When $b=a\Lambda_{j,k}$, the constants below are combined with
Proposition~\ref{prop:Lambda-deriv-bounds:full} and
Proposition~\ref{prop:mult-by-Lambda-preserves-S-anisotropic:full}. In that use, all constants are
finite-seminorm constants at the derivative level under consideration.

The Moyal product is denoted by $a\# b$. Formally,
\begin{equation}\label{eq:moyal-exponential}
a\# b
=
\exp\!\Big(\frac{i}{2}\big(\partial_x\!\cdot\!\partial_\eta-
\partial_\xi\!\cdot\!\partial_y\big)\Big)
\big[a(x,\xi)b(y,\eta)\big]\Big|_{(y,\eta)=(x,\xi)},
\end{equation}
and the corresponding asymptotic expansion is
\begin{equation}\label{eq:moyal-series}
a\# b\sim
\sum_{r=0}^{\infty}\ \sum_{|\alpha|+|\beta|=r}
\frac{(i/2)^r}{\alpha!\,\beta!}
\big(\partial_\xi^\alpha\partial_x^\beta a\big)
\big(\partial_x^\alpha\partial_\xi^\beta b\big).
\end{equation}

\begin{theorem}[Weyl composition in $S^{m_1,m_2}$]\label{thm:weyl-moyal-anisotropic}
Let $a\in S^{m_{1,a},m_{2,a}}$ and $b\in S^{m_{1,b},m_{2,b}}$. For every $N\in\mathbb{N}$,
\[
a\# b
=\sum_{|\alpha|+|\beta|<N}
\frac{(i/2)^{|\alpha|+|\beta|}}{\alpha!\,\beta!}
\big(\partial_\xi^\alpha\partial_x^\beta a\big)
\big(\partial_x^\alpha\partial_\xi^\beta b\big)+r_N,
\]
where
\[
r_N\in S^{m_{1,a}+m_{1,b}-N,\,m_{2,a}+m_{2,b}-N}.
\]
Moreover, for every pair of seminorm indices $(\mu,\nu)$ there exists
$M=M(N,\mu,\nu,n)$ such that
\begin{equation}\label{eq:weyl-remainder-bound}
\|r_N\|_{\mu,\nu}^{(m_{1,a}+m_{1,b}-N,\,m_{2,a}+m_{2,b}-N)}
\leq C_{N,\mu,\nu}
\sum_{|\alpha|+|\beta|\leq M}\|a\|_{\alpha,\beta}^{(m_{1,a},m_{2,a})}
\sum_{|\alpha'|+|\beta'|\leq M}\|b\|_{\alpha',\beta'}^{(m_{1,b},m_{2,b})}.
\end{equation}
Consequently,
\begin{equation}\label{eq:weyl-composition}
\operatorname{Op}^w(a)\operatorname{Op}^w(b)
=\operatorname{Op}^w\!\left(
\sum_{|\alpha|+|\beta|<N}
\frac{(i/2)^{|\alpha|+|\beta|}}{\alpha!\,\beta!}
\big(\partial_\xi^\alpha\partial_x^\beta a\big)
\big(\partial_x^\alpha\partial_\xi^\beta b\big)
\right)+\operatorname{Op}^w(r_N).
\end{equation}
If $\chi,\chi'\in C_c^\infty(K)$ for a compact set $K\subset\mathbb{R}^n$, then
$\chi\operatorname{Op}^w(r_N)\chi'$ has local order $m_{2,a}+m_{2,b}-N$ in the frequency variable.
In particular, by choosing $N$ arbitrarily large, the localized remainders are smoothing.
\end{theorem}

\begin{remark}[Uniformity over mesoscopic patches]\label{rem:uniform-parches}
Applied to $b=a\Lambda_{j,k}$, Theorem~\ref{thm:weyl-moyal-anisotropic} is used only at a prescribed
finite derivative level.  For the mesoscopic partition, Proposition~\ref{prop:mult-by-Lambda-preserves-S-anisotropic:full}
gives, on every fixed compact $K\subset\mathbb{R}^n$, the sharper scale-adapted bound
\[
|\partial_x^\alpha\partial_\xi^\beta(a\Lambda_{j,k})(x,\xi)|
\leq
C_{\alpha,\beta,K}
\mathcal S_{A,B}(a;K)
R_k^{2m_2+|\alpha|-|\beta|}.
\]
Thus one $x$-derivative costs one factor $R_k$, whereas one $\xi$-derivative gains one factor
$R_k^{-1}$.  This finite-level statement is uniform in $(j,k)$ and is the estimate used in the
Sobolev-conjugation and almost-orthogonality arguments.  It should not be replaced by a claim that the
whole family is uniformly bounded, over all derivative orders, in a single fixed global symbol class.
\end{remark}

For the semiclassical calculus we write $S_h^{m_1,m_2}$ for the set of symbols
$a(x,\xi;h)$, $0<h\leq1$, such that for every $\alpha,\beta$
\begin{equation}\label{eq:semi-class}
\sup_{0<h\leq1}\|a(\cdot;h)\|_{\alpha,\beta}^{(m_1,m_2)}<\infty.
\end{equation}
The semiclassical Weyl quantization is
\[
\operatorname{Op}_h^w(a)u(x)=(2\pi h)^{-n}\iint
 e^{\frac{i}{h}(x-y)\cdot\xi}a\!\left(\frac{x+y}{2},\xi;h\right)u(y)\,\mathrm{d}y\,\mathrm{d}\xi.
\]

\begin{theorem}[Semiclassical composition, with controlled remainder]\label{thm:weyl-moyal-semiclass}
Let $a\in S_h^{m_{1,a},m_{2,a}}$ and $b\in S_h^{m_{1,b},m_{2,b}}$. For every $N\in\mathbb{N}$,
\[
a\#_h b
=\sum_{|\alpha|+|\beta|<N}\frac{(ih/2)^{|\alpha|+|\beta|}}{\alpha!\,\beta!}
\big(\partial_\xi^\alpha\partial_x^\beta a\big)
\big(\partial_x^\alpha\partial_\xi^\beta b\big)+h^N r_N(\cdot;h),
\]
where
$r_N(\cdot;h)\in S_h^{m_{1,a}+m_{1,b}-N,\,m_{2,a}+m_{2,b}-N}$ and, for every $(\mu,\nu)$,
\begin{align}
\sup_{0<h\leq1}\|r_N(\cdot;h)\|_{\mu,\nu}^{(m_{1,a}+m_{1,b}-N,\,m_{2,a}+m_{2,b}-N)}
&\leq C_{N,\mu,\nu}
\sum_{|\alpha|+|\beta|\leq M}\sup_{0<h\leq1}\|a(\cdot;h)\|_{\alpha,\beta}^{(m_{1,a},m_{2,a})} \nonumber\\
&\quad\times
\sum_{|\alpha'|+|\beta'|\leq M}\sup_{0<h\leq1}\|b(\cdot;h)\|_{\alpha',\beta'}^{(m_{1,b},m_{2,b})}.
\label{eq:weyl-remainder-semiclass}
\end{align}
Equivalently,
\begin{equation}\label{eq:weyl-composition-semiclass}
\operatorname{Op}_h^w(a)\operatorname{Op}_h^w(b)
=\operatorname{Op}_h^w\!\left(
\sum_{|\alpha|+|\beta|<N}\frac{(ih/2)^{|\alpha|+|\beta|}}{\alpha!\,\beta!}
\big(\partial_\xi^\alpha\partial_x^\beta a\big)
\big(\partial_x^\alpha\partial_\xi^\beta b\big)
\right)+h^N\operatorname{Op}_h^w(r_N(\cdot;h)).
\end{equation}
After inserting compact cutoffs in $x$, the remainder has local frequency order
$m_{2,a}+m_{2,b}-N$, uniformly for $0<h\leq1$.
\end{theorem}

\begin{remark}[Compatibility with the mesoscopic Sobolev reduction]\label{rem:compat-cv-local}
In the body of the paper, the Sobolev-conjugated block is treated at the scale
$R_k=2^{k/2}$.  If $R=R_k$ and
\[
(x,\xi)=(X/R,R\Xi),
\]
then the normalized block symbol has uniformly bounded finite derivatives in $(X,\Xi)$.  The rescaled
spatial cutoffs $\chi(X/R)$ and $\chi'(X/R)$ gain one factor $R^{-1}$ under every nonzero $X$-derivative,
while the normalized Sobolev weights gain one factor $R^{-1}$ under every nonzero $\Xi$-derivative on the
rescaled block, where $|\Xi|\asymp R$.  Consequently every nontrivial Weyl contraction between a block
factor and one of the external cutoff or Sobolev factors gains at least one factor $R^{-1}$.

Lemma~\ref{lem:mesoscopic-kernel-localization} uses the exact oscillatory formula rather than an infinite
formal expansion: a finite number of integrations by parts gives the required phase-space tails and keeps
the constants uniform in $(j,k)$.  The explicit sufficient derivative budget used for the
Cotlar--Stein exponent $M_{\operatorname{CS}}=2n+2$ is recorded below in
Proposition~\ref{prop:mesoscopic-cotlar-derivative-budget}.
\end{remark}

The symbolic composition formulae and the corresponding remainder estimates are standard; see, for example,
\cite[Section~18.5]{HormanderIII}, \cite[Section~14]{TaylorPDEII}, and, for the semiclassical version,
\cite[Section~4.3]{ZworskiSemiclassical}. The adaptation to $S^{m_1,m_2}$ is obtained by the same proof with
the two weights $\langle x\rangle$ and $\langle\xi\rangle$ carried through the seminorm estimates.

\subsection{Cotlar--Stein lemma and Calder\'on--Vaillancourt theorem}

\begin{lemma}[Cotlar--Stein]\label{lem:cotlar-stein}
Let $\{T_j\}_{j\in J}$ be an at most countable family of bounded operators on a Hilbert space
$\mathcal{H}$. Assume that
\[
A:=\sup_{j\in J}\sum_{k\in J}\|T_j^*T_k\|^{1/2}<\infty,
\qquad
B:=\sup_{k\in J}\sum_{j\in J}\|T_jT_k^*\|^{1/2}<\infty.
\]
Then, for every finite subset $F\subset J$,
\[
\left\|\sum_{j\in F}T_j\right\|_{\mathcal{L}(\mathcal{H})}\leq \sqrt{AB}.
\]
In particular, the net of finite partial sums is uniformly bounded and converges in the strong operator
topology to a bounded operator $T$ with $\|T\|\leq\sqrt{AB}$.
\end{lemma}

\begin{corollary}[Tail criterion for norm convergence]\label{cor:cotlar-stein-norm}
In the setting of Lemma~\ref{lem:cotlar-stein}, suppose that $J$ is equipped with an exhaustion
$J_1\subset J_2\subset\cdots$ by finite subsets and set $J^{(N)}:=J\setminus J_N$. If
\[
A_N:=\sup_{j\in J}\sum_{k\in J^{(N)}}\|T_j^*T_k\|^{1/2}\longrightarrow0,
\qquad
B_N:=\sup_{k\in J}\sum_{j\in J^{(N)}}\|T_jT_k^*\|^{1/2}\longrightarrow0,
\]
then $\sum_{j\in J}T_j$ converges in operator norm and the limit satisfies
$\|\sum_{j\in J}T_j\|\leq\sqrt{AB}$.
\end{corollary}

\begin{theorem}[Calder\'on--Vaillancourt, local version]\label{thm:CV-local}
Let $n\geq1$. There exists an integer $M=M(n)$ such that the following holds. Let
$K\subset\mathbb{R}^n$ be compact and let $\chi,\chi'\in C_c^\infty(K)$. If
$b\in S^m_{1,0}(\mathbb{R}^n_x\times\mathbb{R}^n_\xi)$ in the standard sense locally over $K$, namely
\[
B_m:=\max_{|\alpha|+|\beta|\leq M}
\sup_{x\in K,\,\xi\in\mathbb{R}^n}
\langle\xi\rangle^{-m+|\beta|}
\big|\partial_x^\alpha\partial_\xi^\beta b(x,\xi)\big|<\infty,
\]
then, for every $s\in\mathbb{R}$,
\[
\|\chi\operatorname{Op}^w(b)\chi'\|_{H^s\to H^{s-m}}
\leq C\,B_m.
\]
Here $C$ depends on $n$, $s$, $m$, $K$, and finitely many seminorms of $\chi,\chi'$, but not on $b$ except through $B_m$.
\end{theorem}

\begin{proof}
This is the local form of the Calder\'on--Vaillancourt theorem, followed by the standard Sobolev
conjugation argument for symbols of order $m$; see \cite[Theorem~18.6.1]{TaylorPDEII} and
\cite[Chapter~23]{ShubinPDE}. The compact cutoffs restrict the $x$-variable to $K$, so only the displayed
local $x$-seminorms of $b$ enter the constant.
\end{proof}

\begin{theorem}[Semiclassical Calder\'on--Vaillancourt, uniform in $h$]\label{thm:CV-semiclass-uniform}
Let $a(x,\xi;h)\in C^\infty(\mathbb{R}^{2n}\times(0,1])$ and suppose that, for some
$M=M(n)$,
\[
C_M:=\sup_{0<h\leq1}\max_{|\alpha|+|\beta|\leq M}
\sup_{(x,\xi)\in\mathbb{R}^{2n}}
\big|\partial_x^\alpha\partial_\xi^\beta a(x,\xi;h)\big|<\infty.
\]
Then
\[
\|\operatorname{Op}_h^w(a)\|_{L^2\to L^2}\leq C_n C_M,
\qquad 0<h\leq1.
\]
More generally, after Sobolev conjugation the same conclusion holds on $H^s$ provided the corresponding
finite seminorms of the conjugated symbol are uniformly bounded.
\end{theorem}

\begin{proof}
This is the standard semiclassical Calder\'on--Vaillancourt theorem; see \cite[Proposition~4.12]{ZworskiSemiclassical}.  The final Sobolev statement follows by applying the same $L^2$ estimate to the Weyl symbol of
\[
\langle hD\rangle^s
\operatorname{Op}_h^w(a)
\langle hD\rangle^{-s},
\]
or to the non-semiclassical Sobolev weights after the corresponding finite symbol bounds have been verified.
\end{proof}

\begin{corollary}[Local semiclassical version with cutoffs]\label{cor:CV-local-semiclass}
Let $K\subset\mathbb{R}^n$ be compact and let $\chi,\chi'\in C_c^\infty(K)$. If the finite seminorms
required in Theorem~\ref{thm:CV-semiclass-uniform} are bounded uniformly for $0<h\leq1$ on
$K\times\mathbb{R}^n_\xi$, then
\[
\|\chi\operatorname{Op}_h^w(a)\chi'\|_{L^2\to L^2}
\leq C_K C_M,
\qquad 0<h\leq1.
\]
The analogous $H^s$ estimate holds after verifying the corresponding finite seminorm bounds for the
Sobolev-conjugated symbol.
\end{corollary}

\subsection{Finite derivative bookkeeping for global weighted Sobolev bounds}
\label{subsec:appendix-weighted-sobolev-bookkeeping}

We record a finite-seminorm version of the global product calculus needed
for the weighted Sobolev estimate in
Theorem~\ref{thm:global-weighted-sobolev}.  Let
$M_{\operatorname{CV}}=M_{\operatorname{CV}}(n)$ denote the
finite derivative order appearing in
Theorem~\ref{thm:CV-semiclass-uniform} at $h=1$.

\begin{lemma}[Global Weyl-product bookkeeping]
\label{lem:weyl-product-global-bookkeeping}
Let $F\in S^{p,q}$ and $G\in S^{r,t}$, and fix an output derivative
order $N_0\in\mathbb{N}_0$.  Set
\[
\nu(p,q;r,t)
:=
|p|+|q|+|r|+|t|
\]
and
\[
Q(p,q;r,t)
:=
n+1+
\left\lceil
\frac{\nu(p,q;r,t)}{2}
\right\rceil.
\]
Then the seminorms of $F\#G$ in $S^{p+r,q+t}$ involving derivatives
$\partial_x^\alpha\partial_\xi^\beta$ with
$|\alpha|+|\beta|\leq N_0$ are controlled by finitely many seminorms of
$F$ and $G$ involving derivatives of total order at most
\[
N_0+2Q(p,q;r,t).
\]
More precisely, there exists a constant
$C=C(n,N_0,p,q,r,t)$ such that
\[
\max_{|\alpha|+|\beta|\leq N_0}
\|F\#G\|_{\alpha,\beta}^{(p+r,q+t)}
\leq
C
\max_{|\alpha|+|\beta|\leq N_0+2Q}
\|F\|_{\alpha,\beta}^{(p,q)}
\max_{|\alpha|+|\beta|\leq N_0+2Q}
\|G\|_{\alpha,\beta}^{(r,t)},
\]
where $Q=Q(p,q;r,t)$.
\end{lemma}

\begin{proof}
Use the exact oscillatory formula
\[
(F\#G)(Z)
=
\pi^{-2n}\operatorname{Os}\!\iint
e^{-2i\sigma(U,V)}F(Z+U)G(Z+V)
\,\mathrm{d}U\,\mathrm{d}V.
\]
After differentiating in $Z$ up to total order $N_0$, distribute the
output derivatives between the two factors.  Peetre's inequality, applied
separately to the $x$- and $\xi$-weights, shows that after division by the
target weight of order $(p+r,q+t)$ the translated factors contribute at
most a polynomial factor
\[
C\langle U\rangle^{\nu(p,q;r,t)}
\langle V\rangle^{\nu(p,q;r,t)}.
\]
Apply $Q$ times each of the identities
\[
(1-\Delta_U)e^{-2i\sigma(U,V)}
=
(1+4|V|^2)e^{-2i\sigma(U,V)},
\]
\[
(1-\Delta_V)e^{-2i\sigma(U,V)}
=
(1+4|U|^2)e^{-2i\sigma(U,V)}.
\]
The resulting integrand contains the weights
\[
\langle U\rangle^{-2Q+\nu}
\langle V\rangle^{-2Q+\nu}.
\]
By the definition of $Q$,
\[
2Q-\nu>2n,
\]
so both weights are integrable on $\mathbb{R}^{2n}$.  The regularization
uses at most $2Q$ additional derivatives of each input factor.  This gives
the stated finite-seminorm bound.
\end{proof}

For the five factors appearing in the weighted Sobolev conjugation, set
\[
\nu_2
:=
|\sigma-m_1|+|s-m_2|+|m_1|+|m_2|,
\]
\[
\nu_3
:=
|\sigma|+2|s|,
\qquad
\nu_4:=2|\sigma|,
\]
and define
\[
Q_2
:=
n+1+\left\lceil\frac{\nu_2}{2}\right\rceil,
\qquad
Q_3
:=
n+1+\left\lceil\frac{\nu_3}{2}\right\rceil,
\qquad
Q_4
:=
n+1+\left\lceil\frac{\nu_4}{2}\right\rceil.
\]
Finally, define the sufficient weighted Sobolev derivative budget
\begin{equation}\label{eq:weighted-sobolev-budget}
D_{\operatorname{wSob}}
:=
M_{\operatorname{CV}}+2(Q_2+Q_3+Q_4).
\end{equation}

\begin{proposition}[Explicit finite-seminorm budget for weighted Sobolev conjugation]
\label{prop:weighted-sobolev-derivative-budget}
Let $a\in S^{m_1,m_2}$ and let $s,\sigma\in\mathbb{R}$.  Define
\[
c_{s,\sigma}
:=
\langle x\rangle^{\sigma-m_1}
\#
\langle\xi\rangle^{s-m_2}
\# a\#
\langle\xi\rangle^{-s}
\#
\langle x\rangle^{-\sigma}.
\]
Then $c_{s,\sigma}\in S^{0,0}$ and
\[
\max_{|\alpha|+|\beta|\leq M_{\operatorname{CV}}}
\sup_{(x,\xi)\in\mathbb{R}^{2n}}
\big|
\partial_x^\alpha\partial_\xi^\beta
c_{s,\sigma}(x,\xi)
\big|
\leq
C_{n,s,\sigma,m_1,m_2}
\mathcal S_{D_{\operatorname{wSob}},D_{\operatorname{wSob}}}(a).
\]
In particular, the weighted Sobolev conjugation uses no derivative of
$a$ of order larger than $D_{\operatorname{wSob}}$ in either variable.
\end{proposition}

\begin{proof}
Parenthesize the product as
\[
c_{s,\sigma}
=
\Big(
\big(
(\langle x\rangle^{\sigma-m_1}
\#\langle\xi\rangle^{s-m_2})
\# a
\big)
\#\langle\xi\rangle^{-s}
\Big)
\#\langle x\rangle^{-\sigma}.
\]
The successive symbolic orders are
\[
(\sigma-m_1,s-m_2),
\qquad
(\sigma,s),
\qquad
(\sigma,0),
\qquad
(0,0).
\]
The first product contains no factor $a$, so its derivative cost is
absorbed into the constant depending on $s,\sigma,m_1,m_2$.  To control
the final symbol through total order $M_{\operatorname{CV}}$, apply
Lemma~\ref{lem:weyl-product-global-bookkeeping} backwards through the
remaining three products.  The last product requires the preceding symbol
through order $M_{\operatorname{CV}}+2Q_4$; the previous one then requires
order $M_{\operatorname{CV}}+2Q_4+2Q_3$; and the product with $a$ requires
$a$ through order
\[
M_{\operatorname{CV}}+2(Q_2+Q_3+Q_4)
=
D_{\operatorname{wSob}}.
\]
The order identities above show that the exact final symbol belongs to
$S^{0,0}$.  The displayed estimate follows from the same three
applications of the global bookkeeping lemma.
\end{proof}

\begin{remark}
The integer $D_{\operatorname{wSob}}$ is a sufficient finite-seminorm
budget, not an optimal one.  Its only dimensional input is the finite
Calder\'on--Vaillancourt order $M_{\operatorname{CV}}(n)$ already fixed in
Theorem~\ref{thm:CV-semiclass-uniform}; all additional derivative costs are
made explicit in \eqref{eq:weighted-sobolev-budget}.
\end{remark}

\subsection{Finite derivative bookkeeping for mesoscopic almost orthogonality}
\label{subsec:appendix-mesoscopic-bookkeeping}

The proof of Lemma~\ref{lem:mesoscopic-offdiagonal-cotlar} uses only finitely many derivatives.  We
record here one explicit sufficient budget.  The numbers below are deliberately conservative; no
optimality is claimed.

Fix
\begin{equation}\label{eq:CS-fixed-exponents}
M_{\operatorname{CS}}:=2n+2,
\qquad
N_{\operatorname{ker}}:=M_{\operatorname{CS}}+n+2=3n+4,
\qquad
D_{\operatorname{ker}}:=M_{\operatorname{CS}}+N_{\operatorname{ker}}+2=5n+8.
\end{equation}
For the Sobolev weights set
\[
\nu_s
:=
2\left\lceil |s|+|s-m_2|\right\rceil+2
\]
and define
\begin{equation}\label{eq:oscillatory-regularization-order}
Q_{\operatorname{osc}}
:=
n+1+
\left\lceil
\frac{N_{\operatorname{ker}}+\nu_s}{2}
\right\rceil.
\end{equation}
Finally, put
\begin{equation}\label{eq:CS-derivative-budget}
D_{\operatorname{CS}}
:=
2D_{\operatorname{ker}}
+8Q_{\operatorname{osc}}
+M_{\operatorname{CS}}+2.
\end{equation}

\begin{lemma}[Oscillatory regularization with an explicit derivative count]
\label{lem:weyl-product-localized-bookkeeping}
Let $Z,U,V\in\mathbb{R}^{2n}$ and let
\[
(F\# G)(Z)
=
\pi^{-2n}\operatorname{Os}\!\iint
 e^{-2i\sigma(U,V)}F(Z+U)G(Z+V)\,\mathrm{d}U\,\mathrm{d}V.
\]
Fix integers $r,L\geq0$ and a polynomial-growth exponent $\nu\geq0$, and set
\[
Q:=n+1+\left\lceil\frac{L+\nu}{2}\right\rceil.
\]
Assume that $G$ is localized in a closed set $\Omega\subset\mathbb{R}^{2n}$ and that the derivatives
of $F$ and $G$ of total order at most $r+2Q$ satisfy the bounds needed below.  In addition, assume that
$F$ has relative polynomial growth of order at most $\nu$ with respect to $\Omega$, in the sense that
\[
|\partial_Z^\gamma F(Z)|
\leq
C_\gamma
\langle 1+\operatorname{dist}(Z,\Omega)\rangle^\nu,
\qquad
|\gamma|\leq r+2Q.
\]
Then every derivative $\partial_Z^\gamma(F\# G)$ with
$|\gamma|\leq r$ is represented by an absolutely convergent integral after applying $Q$ times each of
the identities
\[
(1-\Delta_U)e^{-2i\sigma(U,V)}
=(1+4|V|^2)e^{-2i\sigma(U,V)},
\]
\[
(1-\Delta_V)e^{-2i\sigma(U,V)}
=(1+4|U|^2)e^{-2i\sigma(U,V)}.
\]
Moreover, if $D=\operatorname{dist}(Z,\Omega)$, the same regularization yields a factor
$\langle D\rangle^{-L}$.  Thus one Weyl product at output derivative order $r$ and decay order $L$
uses no more than $r+2Q$ derivatives of each input factor.
\end{lemma}

\begin{proof}
After applying the displayed identities $Q$ times, the oscillatory integral is a finite sum of terms
with weights
\[
\langle U\rangle^{-2Q}\langle V\rangle^{-2Q}
\]
and derivatives of the input factors of order at most $r+2Q$.  Since $U,V\in\mathbb{R}^{2n}$, the
choice $Q>n$ makes both weights integrable.  If $G(Z+V)\neq0$, then $|V|\geq D$ whenever $Z$ is at
distance $D$ from $\Omega$.  Moreover,
\[
\operatorname{dist}(Z+U,\Omega)
\leq D+|U|,
\]
so Peetre's inequality bounds the relative growth of $F(Z+U)$ by
$C\langle D\rangle^\nu\langle U\rangle^\nu$.  We therefore extract $L+\nu$ powers from the
$V$-weight:
\[
\langle V\rangle^{-2Q}
\leq
C_{L,\nu}\langle D\rangle^{-L-\nu}
\langle V\rangle^{-2Q+L+\nu}.
\]
The factor $\langle D\rangle^{-\nu}$ cancels the possible relative growth of $F$, while
$\langle U\rangle^\nu$ is absorbed by the $U$-weight.  By the definition of $Q$,
\[
2Q-L-\nu>2n,
\qquad
2Q-\nu>2n,
\]
so both remaining weights are integrable.  Differentiating in $Z$ before regularization only
adds the prescribed $r$ derivatives.  This proves the claim.
\end{proof}

\begin{proposition}[Explicit finite-seminorm budget for the Cotlar--Stein proof]
\label{prop:mesoscopic-cotlar-derivative-budget}
Let $a\in S^{m_1,m_2}$, fix $s\in\mathbb{R}$ and fixed compactly supported spatial cutoffs as in
Definition~\ref{def:sobolev-conjugated-blocks}.  Let $K_1$ be a fixed compact neighborhood containing
their supports and all midpoint sets occurring in the Weyl kernels.  With the integers in
\eqref{eq:CS-fixed-exponents}--\eqref{eq:CS-derivative-budget}, the proof of
Lemma~\ref{lem:mesoscopic-kernel-localization} and
Lemma~\ref{lem:mesoscopic-offdiagonal-cotlar} at the exponent
$M=M_{\operatorname{CS}}$ is controlled by
\[
\mathcal S_{D_{\operatorname{CS}},D_{\operatorname{CS}}}(a;K_1)
\]
and by
\[
\max_{|\gamma|\leq D_{\operatorname{CS}}}
\sup_{x\in K_1}
\left(
\|\partial_x^\gamma T_x\|_{\operatorname{op}}
+
\|\partial_x^\gamma T_x^{-1}\|_{\operatorname{op}}
\right),
\]
together with finitely many seminorms of the fixed cutoffs up to order
$D_{\operatorname{CS}}$.  In particular, no derivative of $a$, $T$, $T^{-1}$, or the fixed cutoffs of
order larger than $D_{\operatorname{CS}}$ is needed for the quantitative Cotlar--Stein conclusion with
$M_{\operatorname{CS}}=2n+2$.

More precisely, one may organize the proof as follows:
\begin{enumerate}[(i)]
\item use the exact Weyl formula four times for
\[
w_+\#\psi\#(a\Lambda_{j,k})\#\psi'\#w_-,
\]
regularizing each product with at most $Q_{\operatorname{osc}}$ applications of each second-order
operator in Lemma~\ref{lem:weyl-product-localized-bookkeeping};
\item request derivatives of the resulting conjugated symbol through order
$D_{\operatorname{ker}}$ in each variable and phase-space decay of order $N_{\operatorname{ker}}$;
\item obtain the kernel envelope and its $z$-derivatives through order $M_{\operatorname{CS}}$ by at most
$N_{\operatorname{ker}}$ integrations by parts in the centered frequency variable;
\item perform exactly $M_{\operatorname{CS}}=2n+2$ integrations by parts in the interaction variable
$z$ in the nonstationary region.
\end{enumerate}
After extracting the factor
$\operatorname{dist}_{\mathrm{mes}}(\lambda,\mu)^{-M_{\operatorname{CS}}}$, the choice
$N_{\operatorname{ker}}=M_{\operatorname{CS}}+n+2$ leaves a kernel envelope of order $n+2>n$,
which is integrable and therefore sufficient for Schur's test.  Since
$M_{\operatorname{CS}}/2=n+1>n$, the discrete packing sum is Cotlar--Stein summable.
\end{proposition}

\begin{proof}
For the exact Sobolev-conjugated symbol, work in the rescaled variables
\[
X=R_kx,
\qquad
\Xi=R_k^{-1}\xi.
\]
After removal of the scalar order, the localized factor $a\Lambda_{j,k}$ has uniformly bounded
finite derivatives.  The rescaled spatial cutoffs are slowly varying, and the two normalized Sobolev
weights have polynomial growth controlled by $\nu_s$.  Apply
Lemma~\ref{lem:weyl-product-localized-bookkeeping} successively to the four Weyl products.  At output
order at most $2D_{\operatorname{ker}}$ in the combined phase-space variables, each product increases
the required input order by at most $2Q_{\operatorname{osc}}$.  Four products therefore require no more
than
\[
2D_{\operatorname{ker}}+8Q_{\operatorname{osc}}
\]
derivatives of the original factors.  The additional term
$M_{\operatorname{CS}}+2$ in \eqref{eq:CS-derivative-budget} leaves room for the derivatives of the
moving center $T_x^{-1}\zeta_{j,k}$ and for the derivatives of the reciprocal phase gradient occurring
in the $z$-integration by parts.  Thus the displayed $D_{\operatorname{CS}}$ is a sufficient common
upper bound for all derivatives entering the argument.

To obtain the kernel decay, use the centered variable
\[
\xi=T_m^{-1}\zeta_{j,k}+R_k\eta.
\]
The symbol estimate through order $D_{\operatorname{ker}}$ and decay order
$N_{\operatorname{ker}}$ allows the integrations by parts in $\eta$ needed for the kernel envelope and
its derivatives through order $M_{\operatorname{CS}}$.  In the complement of the effective interaction
region, extracting $M_{\operatorname{CS}}$ powers of the normalized block distance leaves the exponent
\[
N_{\operatorname{ker}}-M_{\operatorname{CS}}=n+2,
\]
which is integrable.  In the effective interaction region, the lower bound on the directional phase
gradient from Lemma~\ref{lem:mesoscopic-offdiagonal-cotlar} permits exactly
$M_{\operatorname{CS}}$ integrations by parts.  The resulting Schur estimate is
\[
\|\widetilde T_\lambda^*\widetilde T_\mu\|_{L^2\to L^2}
+
\|\widetilde T_\lambda\widetilde T_\mu^*\|_{L^2\to L^2}
\leq
C
\left(1+\operatorname{dist}_{\mathrm{mes}}(\lambda,\mu)\right)^{-M_{\operatorname{CS}}}.
\]
Finally, the packing estimate in Section~\ref{subsec:dyadic-decomposition} gives summability because
$M_{\operatorname{CS}}/2=n+1>n$.  This proves the claimed finite-seminorm budget.
\end{proof}

\begin{remark}[Role of the derivative budget]
The integer $D_{\operatorname{CS}}$ is a sufficient bookkeeping bound, not an optimal one.  Its purpose
is to make explicit that the off-diagonal estimate and the Cotlar--Stein summability in
Proposition~\ref{thm:global-sum-semiclassical} use only finitely many seminorms, with an order depending
only on $n,s,m_2$ and the fixed cutoffs, and not on the block indices $j,k$.  In particular, the
almost-orthogonality argument does not rely on an implicit all-orders uniform symbol bound.
\end{remark}